\documentclass{amsart}

\usepackage{dha}

\renewcommand{\gamma}{\omega}

\author{Fernando Muro}
\email{fmuro@us.es}
\address{Universidad de Sevilla,
Facultad de Matemáticas,
Departamento de Álgebra,
Avda. Reina Mercedes s/n,
41012 Sevilla, Spain}

\title{Derived universal Massey products}

\begin{document}

\begin{abstract}
    We define an obstruction to the formality of a differential graded algebra over a graded operad defined over a commutative ground ring. This obstruction lives in the derived operadic cohomology of the algebra. Moreover, it determines all operadic Massey products induced on the homology algebra, hence the name of derived universal Massey product.
\end{abstract}

\subjclass[2020]{18M70,55S20,13D03,16E40,17B56}

\keywords{Operad, algebra, cohomology, formality, minimal model, resolution}




\thanks{The author was partially supported by the Spanish Ministry of Economy under the grant MTM2016-76453-C2-1-P (AEI/FEDER, UE) and by the Andalusian Ministry of Economy and Knowledge and the Operational Program FEDER 2014–2020 under the grant US-1263032. The author is very grateful to Steffen Sagave for his useful comments and remarks on a preliminary version of this paper.}

\maketitle

\tableofcontents

\section*{Introduction}

Massey products are a well-known secondary operation on the homology $H_*(A)$ of a differential graded associative $\kk$-algebra $A$ which is defined for any three homogeneous elements $x,y,z\in H_*(A)$ such that $x\cdot y=0=y\cdot z$. It yields an element of the quotient
\[\langle x,y,z\rangle\in \frac{H_{\abs{x}+\abs{y}+\abs{z}+1}(A)}{x\cdot H_{\abs{y}+\abs{z}+1}(A)+H_{\abs{x}+\abs{y}+1}(A)\cdot z},\]
often regarded as a subset (coset) of the numerator, and the denominator is called the indeterminacy.

When $\kk$ is a field, the universal Massey product of $A$ is the Hochschild cohomology class
\[\g{A}{}=\{m_3\}\in HH^{3,-1}(H_*(A),H_*(A))\]
of the first non-trivial piece of a minimal model for $A$, which is a minimal $A$-infinity algebra
\[(H_*(A),m_3,m_4,\dots,m_n,\dots)\]
built upon the homology algebra $H_*(A)$ which is $\infty$-quasi-isomorphic to $A$. Here $m_3$ is an always-defined fully-determined operation
\[m_3\colon H_*(A)\otimes H_*(A)\otimes H_*(A)\To H_*(A)\]
of degree $1$ which determines all Massey products since
\[m_3(x,y,z)\in \langle x,y,z\rangle\]
whenever the Massey product is defined. Despite the choice, the universal Massey product is a quasi-isomorphism invariant of $A$. Moreover, it is an obstruction to formality. This class goes back to Kadeishvili's \cite{kadeishvili_1988_structure_infty_algebra}. It was studied in depth by Benson, Krause, and Schwede in \cite{benson_krause_schwede_2004_realizability_modules_tate}. It has also been considered by Kaledin in \cite{kaledin_2007_remarks_formality_families}. The name we give it here is inspired by the similar notion of universal Toda bracket due to Baues and his collaborators \cite{baues_dreckmann_1989_cohomology_homotopy_categories,baues_tonks_1996_sumnormalised_cohomology_categories}.

Dimitrova \cite{dimitrova_2012_obstruction_theory_operadic} extended the definition of the universal Massey product to differential graded algebras $A$ over a graded operad $\O{O}$ defined over a ground field $\kk$. Her definition simplifies if $\O{O}$ is a Koszul operad in characteristic zero or a non-symmetric Koszul operad in positive characteristic. In that case, the universal Massey product is an invariant living in operadic cohomology in the sense of \cite{loday_vallette_2012_algebraic_operads},
\[\g{A}{\O{O}}=\{\m\}\in H^{2,-1}_{\O{O}}(H_*(A),H_*(A)),\]
it is an obstruction to formality,
the representing cocycle $\m$ is defined from a minimal $\Oinf{\O{O}}$-model of $A$, and any representing cocycle determines all operadic Massey products \cite{muro_2021_massey_products_algebras}. This generalizes the associative case (operadic cohomology over the associative operad essentially coincides with Hochschild's up to a degree shift).

In this paper we extend the theory of universal Massey products to algebras $A$ over a graded operad $\O{O}$ over a commutative ground ring $\kk$. We assume that $\O{O}$ is a Koszul operad and $\Q\subset \kk$ or that $\O{O}$ is a non-symmetric Koszul operad with no restrictions on $\kk$. Under these standing assumptions, we consider derived operadic cohomology. Here the word ``derived" has the usual homotopical meaning. Indeed, the category of differential graded $\O{O}$-algebras carries a model structure with quasi-isomorphisms as weak equivalences and surjections as fibrations. Hence derived operadic cohomology can be constructed by means of cofibrant replacements. We define an invariant derived operadic cohomology class in the expected bidegree
\[\g{A}{\O{O}}=\{\m\}\in \D{H}^{2,-1}_{\O{O}}(H_*(A),H_*(A)).\]
This class coincides with Dimitrova's for $\kk$ a field (there is no need to derive in this case) and, in general, it computes the operadic Massey products defined in \cite{muro_2021_massey_products_algebras}, so it deserves to be called \emph{derived universal Massey product}. As in previous cases, it is also an obstruction to formality.

A representing cocycle $\m$ is defined from a derived minimal model in the sense of \cite{derived_homotopy_algebras}. These minimal models extend the classical theory over fields to commutative rings. They were first defined by Sagave \cite{sagave_2010_dgalgebras_derived_algebras} in the associative case. He also defined a derived universal Massey product living in a cohomology theory which is unfortunately not invariant (it is \emph{dependent} on the choice of resolutions used to define it, see Remark \ref{error_sagave}). Nevertheless, this cohomology is not that far away from derived Hochschild cohomology (also known as Shukla cohomology), and our $\m$ in the associative case is a slight variation (actually a simplification) of Sagave's cocycle (see Remark \ref{classical_derived_universal_massey_product}). Therefore, this paper can be regarded as a generalization of Sagave's derived universal Massey product to algebras over Koszul operads after some fixing.


We start this paper with a few elementary facts on derived operadic cohomology (\S\ref{section:derived_operadic_cohomology}). We then split this cohomology for graded algebras (\S\ref{section:splitting_cohomology}) using certain bicomplex-shaped resolutions arising from the Cartan--Eilenberg homotopy theory of algebras in bicomplexes \cite[\S2]{derived_homotopy_algebras}. Derived universal Massey products will later be defined in one of these pieces (\S\ref{section:universal_massey_products}). In order to achieve this goal, we recall in \S\ref{section:derived_homotopy_algebras} what we need from the theory of derived minimal models.

In \S\ref{section:chain_complexes} we consider the case where $\O{O}$ is the initial operad $\unit$. Algebras over the initial operad are just graded modules (or complexes if they are equipped with a differential). This case is interesting because, by invariance, it sometimes allows to show that the universal Massey product of a differential graded $\O{O}$-algebra $A$ (with $\O{O}$ any other admissible operad) is non-trivial by just looking at its underlying complex. More precisely, it suffices that some
\[H_{n+1}(A)\hookrightarrow\frac{A_{n+1}}{d(A_{n+2})}\stackrel{d}{\To}\ker[A_{n}\stackrel{d}{\to} A_{n-1}]\twoheadrightarrow H_n(A)\]
is a non-trivial $\kk$-module extension. However, this can only happen if $\kk$ has global dimension $\geq 2$.

Finally, in \S\ref{section:torsion_massey}, we show examples of $\O{O}$-algebras with $\kk$ hereditary (i.e.~of global dimension $\leq 1$) which have non-trivial derived universal Massey product. For this, we establish a connection between  derived universal Massey products and certain operadic Massey products, in the sense of \cite{muro_2021_massey_products_algebras}, associated to torsion elements.


Throughout this paper, $\kk$ is a commutative ground ring and $\O{O}=\P{E}{R}$ is a graded $\kk$-projective quadratic Koszul operad with generating $\SS$-module $E$, which is trivial in arity $0$, i.e.~$E(0)=0$, and relations $\SS$-module $R\subset E\circ_{(1)}E$. We refer the reader to \cite{fresse_2004_koszul_duality_operads} for operadic Koszul duality over commutative rings. Nonetheless, we use terminology and notation from \cite{loday_vallette_2012_algebraic_operads} since we feel it is more commonly employed.
We stress that $\O{O}$ is \emph{graded}, i.e.~it has trivial differential, hence differential graded $\O{O}$-algebras $A$ have graded homology $\O{O}$-algebras $H_*(A)$. The operad $\O{O}$ is naturally weighted and we denote by $\O{O}^{(n)}$ its weight $n$ part. The weighting is a non-negative grading characterized by the fact that $E$ has weight $1$. We have $\O{O}^{(0)}=\unit$, the $\SS$-module given by $\kk$ concentrated in arity $1$ and degree $0$, and $1\in \kk=\unit(1)_0\subset\O{O}(1)_0$ is the operadic unit. The  $\SS$-module $\unit$ is the monoidal unit for the composition product whose (co)monoids are (co)operads, hence it is the initial operad (and the final cooperad). We also have $\O{O}^{(1)}=E$ and $\O{O}^{(2)}=(E\circ_{(1)}E)/R$. We assume that either $\Q\subset\kk$ or $\O{O}$ is non-symmetric (or rather the symmetrization of a non-symmetric operad). In this way we can apply all results from \cite{derived_homotopy_algebras}.

\section{Derived operadic cohomology}\label{section:derived_operadic_cohomology}

In this section we recall operadic cohomology from \cite{loday_vallette_2012_algebraic_operads} and consider its derived functor. By \cite{hinich_1997_homological_algebra_homotopy,hinich_2003_erratum_homological_algebra,muro_2011_homotopy_theory_nonsymmetric,muro_2017_correction_articles_homotopy}, under our standing assumptions, the category of differential graded $\O{O}$-algebras carries a model structure transferred from the category $\chain$ of chain complexes. Fibrations are surjections and weak equivalences are quasi-isomorphisms, respectively.

The \emph{suspension} $\suspension X$ of a chain complex $X$ is defined by shifting one degree up and changing the sign of the differential. Equivalently, $\suspension$ is defined as a functor by the existence of a natural degree $1$ isomorphism $\suspension\colon X\to \suspension X$.

The \emph{Koszul dual cooperad} $\K{\O{O}}=\C{sE}{s^2R}$ is cogenerated by $sE$ with corelations $s^2R\subset s^2(E\circ_{(1)}E)\cong(sE)\circ_{(1)}(sE)$. It is naturally weighted with $(\K{\O{O}})^{(0)}=\unit$, $(\K{\O{O}})^{(1)}=sE$, and $(\K{\O{O}})^{(2)}=s^2R$. We denote by
\[\kappa\colon\K{\O{O}}\To\O{O}\]
the composite
\[\K{\O{O}}\twoheadrightarrow sE\stackrel{s^{-1}}{\To}E\hookrightarrow\O{O},\]
where the first arrow is the projection onto the weight $1$ part of $\K{\O{O}}$, the second one is the desuspension, and the third one is the inclusion of the weight $1$ part of $\O{O}$.

Given $\mu\in\K{\O{O}}(r)$ we use the following Sweedler-like formula for its infinitesimal decomposition, like in \cite[\S10.1.2]{loday_vallette_2012_algebraic_operads},
\begin{equation}\label{sweedler_infinitesimal_decomposition}
    \Delta_{(1)}\colon\K{\O{O}}\To \K{\O{O}}\circ_{(1)}\K{\O{O}},\qquad
    \Delta_{(1)}(\mu)=\sum_{(\mu)}(\mu^{(1)}\circ_{l}\mu^{(2)})\cdot\sigma.
\end{equation}
Here $\sigma\in\SS_r$ is a permutation.

\begin{definition}[{\cite[\S12.4]{loday_vallette_2012_algebraic_operads}}]\label{operadic_cochain_complex}
    The \emph{operadic cochain complex} $C_{\O{O}}^*(A,M)$ of a differential graded $\O{O}$-algebra $A$ with coefficients in an $A$-module $M$ is given by
    \[C_{\O{O}}^n(A,M)=\prod_{p\in\Z}\hom(\K{\O{O}}(A)_p, M_{p-n}).\]
    The degree $+1$ differential is defined as follows for $\psi\in C_{\O{O}}^n(A,M)$,
    \begin{multline}\label{operadic_cohomology_differential_0}
        d(\psi)(\mu;x_1,\dots,x_r) ={} d(\psi(\mu;x_1,\dots,x_r)                         )                                                                                   \\
        +\sum_{(\mu)}(-1)^{\eta_n}\kappa(\mu^{(1)})(x_{\sigma^{-1}(1)},\dots,\psi(\mu^{(2)};x_{\sigma^{-1}(l)},\dots),\dots)               \\
        -\sum_{(\mu)}(-1)^{n+\eta_1}\psi(\mu^{(1)};x_{\sigma^{-1}(1)},\dots,\kappa(\mu^{(2)})(x_{\sigma^{-1}(l)},\dots),\dots) \\
        -\sum_{i=1}^r(-1)^{n+\beta}\psi(\mu;x_1,\dots,d(x_i),\dots,x_r).
    \end{multline}
    Here
    \begin{equation}\label{constants_beta_eta_alpha}
        \begin{split}
            \eta_n & =\alpha_{\sigma}+\abs{\mu^{(1)}}n+(\abs{\mu^{(2)}}+n)\sum_{m=1}^{l-1}\abs{x_{\sigma^{-1}(m)}},\qquad \beta=\abs{\mu}+\sum_{j=1}^{i-1}\abs{x_j},\\
            \alpha_\sigma  &=\sum_{\substack{s<t                                                      \\\sigma(s)>\sigma(t)}}\abs{x_s}\abs{x_t}.
        \end{split}
    \end{equation}

    The cohomology of $C_{\O{O}}^*(A,M)$ is called \emph{operadic cohomology} and denoted by
    \[H_{\O{O}}^*(A,M).\]
    The \emph{derived operadic cohomology} of $A$ with coefficients in $M$ is defined as
    \[\D{H}^{*}_{\O{O}}(A,M)=H^{*}_{\O{O}}(B,M),\]
    where $B$ is a cofibrant resolution of $A$, given by a trivial fibration of $\O{O}$-algebras $B\to A$ with cofibrant source which turns $M$ into a $B$-module. This trivial fibration induces a ``non-derived to derived" comparison map
    \[H^{*}_{\O{O}}(A,M)\To \D{H}^{*}_{\O{O}}(A,M).\]

    The bivariant functoriality of (derived) operadic cohomology is clear and the comparison map is natural.
\end{definition}

\begin{remark}\label{cohomology_filtration}
    The operadic cochain complex admits a decreasing filtration by the weight in $\K{\O{O}}$, i.e.~$\psi\in F_nC_{\O{O}}(A,M)$ if $\psi(\mu;x_1,\dots,x_r)=0$ for $\mu$ of weight $<n$,
    \[\cdots\subset F_{n+1}C_{\O{O}}(A,M)\subset F_nC_{\O{O}}(A,M)\subset\cdots.\]
    Since $\kappa$ vanishes except in weight $1$, $\mu^{(2)}$ and $\mu^{(1)}$ in the second and third lines of \eqref{operadic_cohomology_differential_0} have strictly less weight than $\mu$. This shows that, if $\psi\in F_n$, the second and third lines of \eqref{operadic_cohomology_differential_0}, are cochains in $F_{n+1}$, and the first and last lines are cochains in $F_n$. Therefore the previous filtration is indeed a filtration and, moreover,
    \[E_0(C_{\O{O}}(A,M))=\hom_{\chain}(\K{\O{O}}(A),M),\]
    the inner hom in $\chain$, which does not take into account either the algebra structure of $A$ or the module structure of $M$. The filtration is exhaustive since $F_0=C_{\O{O}}(A,M)$. It is not bounded above, but it is weakly convergent in the sense of \cite[Definition 3.1]{mccleary_2001_user_guide_spectral} since $\bigcap_{n}F_{n}C_{\D{\O{O}}}(A,M)=0$. In addition, it is complete in the sense of \cite[Definition 3.8]{mccleary_2001_user_guide_spectral} by the definition of $C_{\O{O}}(A,M)$ as a product, see \cite[Lemma 3.10]{mccleary_2001_user_guide_spectral}. Hence, we can apply the following useful lemma to it.
\end{remark}

\begin{lemma}[{\cite[Theorem 3.9]{mccleary_2001_user_guide_spectral}}]\label{spectral}
    If $f\colon X\to Y$ is a morphism of complete and exhaustive filtered complexes which induces an isomorphism on some page $E^n$ of the corresponding spectral sequences, then $f$ is a quasi-isomorphism.
\end{lemma}

\begin{lemma}\label{cohomology_quasi-iso}
    Given a quasi-isomorphism $f\colon B\to A$ between differential graded $\O{O}$-algebras with underlying cofibrant complexes and a $A$-module $M$, the induced morphism in cohomology $f^*\colon H^*_{\O{O}}(A,M)\to H^*_{\O{O}}(B,M)$ is an isomorphism.
\end{lemma}

\begin{proof}
    The induced map on complexes $f^*\colon C^*_{\O{O}}(A,M)\to C^*_{\O{O}}(B,M)$ is obviously compatible with the filtration in Remark \ref{cohomology_filtration}. The induced morphism between $E_0$ terms $\hom_{\chain}(\K{\O{O}}(A),M)\to \hom_{\chain}(\K{\O{O}}(B),M)$ is $\hom_{\chain}(\K{\O{O}}(f),M)$. By Lemma \ref{spectral}, it suffices to prove that this morphism is a quasi-isomorphism. Since all objects in $\chain$ are fibrant, it is enough to show that $\K{\O{O}}(f)\colon \K{\O{O}}(B)\to \K{\O{O}}(A)$ is a quasi-isomorphism between cofibrant complexes. This will follow from \cite[Proposition 11.5.3]{fresse_2009_modules_operads_functors} since $\O{O}$ is cofibrant as an $\SS$-module in the sense of \cite[\S11.4]{fresse_2009_modules_operads_functors}. This is indeed a consequence of our standing assumptions since $\O{O}$ is $\kk$-projective and moreover, each $\O{O}(r)$ is $\SS_r$-projective because $\O{O}$ is the symmetrization of a non-symmetric operad or the ground ring contains $\Q$ and Maschke's theorem applies. Any arity-wise $\SS_r$-projective $\SS$-module with trivial differential is cofibrant.
\end{proof}

The same argument shows homotopy invariance in the second variable $M$, but we will not use it.

\begin{corollary}\label{corolario_todo}
    The derived cohomology $\D{H}^{n}_{\O{O}}(A,M)$ of an $\O{O}$-algebra $A$ is independent of choices. Moreover, an $\O{O}$-algebra quasi-isomorphism $B\to A$ induces an isomorphism in derived cohomology $f^*\colon \D{H}^{*}_{\O{O}}(A,M)\cong \D{H}^{*}_{\O{O}}(B,M)$ for any $A$-module $M$. Furthermore, if an $\O{O}$-algebra $B$ has underlying cofibrant complex then the comparison map is an isomorphism $H^{*}_{\O{O}}(B,M)\cong\D{H}^{*}_{\O{O}}(B,M)$. In particular, combining both, for any $A$, a quasi-isomorphism $B\to A$ whose source has underlying cofibrant complex induces an isomorphism $\D{H}^{*}_{\O{O}}(A,M)\cong H^{*}_{\O{O}}(B,M)$.
\end{corollary}

\section{Cohomology of graded algebras}\label{section:splitting_cohomology}

In this section we concentrate in the derived operadic cohomology of \emph{graded} $\O{O}$-algebras and modules, i.e.~with trivial differential. We show how to compute it from bicomplex-shaped resolutions and then we split it by introducing an extra grading.

We consider \emph{bicomplexes} $X$ concentrated in the right half plane, i.e.~$X_{p,q}=0$ for $p<0$, with anti-commuting \emph{horizontal} and \emph{vertical differentials} $d_1$ and $d_0$, respectively,
\[d_1\colon X_{p,q}\To X_{p-1,q},\quad d_0\colon X_{p,q}\To X_{p,q-1},\quad d_1^2=0=d_0^2,\quad d_0d_1+d_1d_0=0.\]
They look like
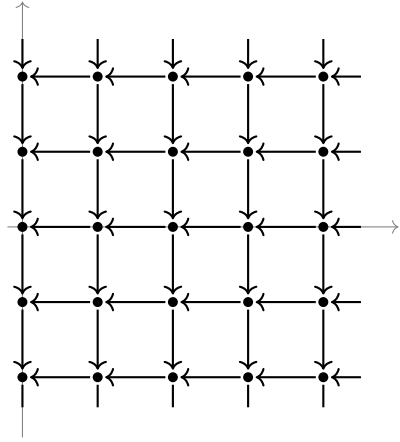
\begin{center}
    \begin{tikzpicture}
        \draw[->, gray] (-0.2,0) -- (5,0); 
        \draw[->, gray] (0,-2.8) -- (0,3); 
        \foreach \x in {0,...,4}
        \foreach \y in {-2,...,2}
            {
                \node[shape=circle,fill=black,scale=.4] (\x-\y) at (\x,\y) {}; 
            }
        \foreach \x in {1,...,4}
        \foreach \y in {-2,...,2}
            {
                \draw[->, thick] (\x-.1,\y) -- (\x-.9,\y); 
            }
        \foreach \y in {-2,...,2}
            {
                \draw[->, thick] (4.5,\y) -- (4.1,\y); 
            }
        \foreach \x in {0,...,4}
        \foreach \y in {-1,...,2}
            {
                \draw[->, thick] (\x,\y-.1) -- (\x,\y-.9); 
            }
        \foreach \x in {0,...,4}
            {
                \draw[->, thick] (\x,2.5) -- (\x,2.1); 
                \draw[thick] (\x,-2.1) -- (\x,-2.4); 
            }
    \end{tikzpicture}
\end{center}
We will not distinguish between a bicomplex and its total complex, i.e.~we regard a bicomplex $X$ as a complex equipped with a splitting,
\begin{equation}\label{bigraded_splitting}
    X_n=\bigoplus_{n=p+q}X_{p,q},
\end{equation}
concentrated in $p\geq 0$, whose differential $d$ decomposes as $d=d_0+d_1$ with $d_0$ and $d_1$ like above. Here $p$ is called \emph{horizontal degree}, $q$ is the \emph{vertical degree}, the pair $(p,q)$ is the \emph{bidegree}, and $n=p+q$ is the \emph{total degree}. \emph{Morphisms} of bicomplexes are chain maps preserving bidegrees. An \emph{$E^2$-equivalence} is a morphism which induces an isomorphism between the $E^2$ pages of the spectral sequences associated to the increasing filtration by the horizontal degree (the first spectral sequence in the sense of \cite[Theorem 2.15]{mccleary_2001_user_guide_spectral}) of the source and target bicomplexes. The category $\bichain$ of bicomplexes is closed symmetric monoidal with the usual tensor product of complexes and the obvious assignation of horizontal and vertical degrees, compare \cite[Definition 2.7]{muro_roitzheim_2019_homotopy_theory_bicomplexes}. The symmetry constraint uses the Koszul sign rule with respect to the total degree. We can regard complexes as bicomplexes concentrated in horizontal degree $0$. This is compatible with the monoidal structures. In this way, we can consider $\O{O}$-algebras in bicomplexes, that we call \emph{bicomplex $\O{O}$-algebras}, see \cite[\S3]{derived_homotopy_algebras}. A bicomplex $\O{O}$-algebra $A$ is just a bicomplex equipped with an $\O{O}$-algebra structure whose structure maps decompose as
\[\O{O}(r)_{s}\otimes A_{p_1,q_1}\otimes\cdots\otimes A_{p_r,q_r} \To A_{p_1+\cdots+p_r,s+q_1+\cdots+q_r}.\]


\begin{definition}\label{horizontal_resolution}
    A bicomplex $\O{O}$-algebra is \emph{minimal} if it has trivial vertical differential, $d_0=0$. A \emph{horizontal resolution} of a graded $\O{O}$-algebra $A$ consists of a $\kk$-projective minimal bicomplex $\O{O}$-algebra $\bires{A}$ equipped with an $E^2$-equivalence $\rho\colon\bires{A}\to A$.
\end{definition}

\begin{remark}\label{certainty}
    In this case, being an $E^2$-equivalence is the same as saying that $\rho$ exhibits $A$ as the horizontal homology of $\bires{A}$. Equivalently, each horizontal complex $(\bires{A})_{*,q}$ is a resolution of $A_q$ via $\rho$.
\end{remark}

The underlying complex of a minimal bicomplex $\O{O}$-algebra is a direct sum of bounded below $\kk$-projective complexes. Such complexes are cofibrant. Therefore, the following result is a direct consequence of Corollary \ref{corolario_todo}.

\begin{corollary}
    If $A$ is a graded $\O{O}$-algebra, $M$ is a graded $A$-module, and $\rho\colon\bires{A}\to A$ is a horizontal resolution then this map induces an isomorphism
    $\D{H}^{n}_{\O{O}}(A,M)\cong H^{n}_{\O{O}}(\bires{A},M)$.
\end{corollary}

We now establish the existence of horizontal resolutions, even a preferred one.

\begin{proposition}\label{horizontal_resolutions_exist}
    Any graded $\O{O}$-algebra $A$ has a horizontal resolution $\rho\colon \bires{A}'\to A$ such that, if $\rho\colon \bires{A}\to A$ is another horizontal resolution, then there exists an $E^2$-equivalence $g\colon \bires{A}'\to \bires{A}$ such that $\rho g=\rho'$. Moreover, if $h\colon A\to B$ is a graded $\O{O}$-algebra morphism and $\rho\colon\bires{B}\to B$ is a horizontal resolution then there exists a bicomplex $\O{O}$-algebra morphism $g\colon \bires{A}'\to \bires{B}$ such that $\rho g=h\rho'$.
\end{proposition}

\begin{proof}
    The category of bicomplex $\O{O}$-algebras can be endowed with the so-called Cartan--Eilenberg semi-model structure \cite[\S2]{derived_homotopy_algebras}.
    A cofibrant replacement $\tilde{A}\to A$ in this semi-model category is an $E^2$-equivalence. The induced morphism between $E^1$ terms is another $E^2$-equivalence $\rho'\colon\bires{A}'=H_*^v(\tilde{A})\to A$. The source is the vertical homology of $\tilde{A}$, which is a minimal bicomplex $\O{O}$-algebra by construction. Moreover, it is $\kk$-projective by \cite[Proposition 2.4]{derived_homotopy_algebras} and \cite[Theorem 4.1]{muro_roitzheim_2019_homotopy_theory_bicomplexes}. Hence $\rho'$ is a horizontal resolution.

    If $\rho$ is another horizontal resolution then $\rho$ is a trivial fibration in the Cartan--Eilenberg semi-model structure by \cite[Theorem 4.1]{muro_roitzheim_2019_homotopy_theory_bicomplexes} and Remark \ref{certainty}. Therefore, there is a lift
    \begin{center}
        \begin{tikzcd}
            &\bires{A}\ar[d,"\rho"]\\
            \tilde{A}\ar[r]\ar[ru,dashed,"\tilde{g}"]&A
        \end{tikzcd}
    \end{center}
    and we can take $g$ to be the map induced by $\tilde{g}$ in vertical homology. The final part of the statement is a slight generalization of this.
\end{proof}

Using horizontal resolutions, we can split the derived operadic cohomology of a graded $\O{O}$-algebra with coefficients in a graded module. In the non-derived case this splitting is straightforward.

\begin{proposition}\label{splitting}
    Given a graded $\O{O}$-algebra $A$, a graded $A$-module $M$, and a horizontal resolution $\bires{A}\to A$, the operadic cochain complex $C_{\O{O}}^*(\bires{A},M)$ splits as follows
    \[C^*_{\O{O}}(\bires{A},M)=\prod_{t\in\Z}s^{-t}C^{*,t}_{\O{O}}(\bires{A},M)\]
    Here $C^{*,t}_{\O{O}}(\bires{A},M)$ is the cochain complex defined as
    \begin{align*}
        C^{w,t}_{\O{O}}(\bires{A},M) 
         & =\bigoplus_{u=0}^w\prod_{p\in\Z}\hom((\K{\O{O}})^{(u)}(\bires{A})_{w-u,p-w+u}, M_{p-w-t}).
    \end{align*}
    Here $(\K{\O{O}})^{(u)}(\bires{A})$ is bigraded because $\bires{A}$ is and we regard $\K{\O{O}}$ as concentrated in horizontal degree $0$.
    The differential of $\psi=(\psi_0,\dots,\psi_w)\in C^{w,t}_{\O{O}}(\bires{A},M)$ is given by
    \begin{multline}\label{operadic_cohomology_differential}
        d(\psi)_u(\mu;x_1,\dots,x_r) = \\
        \sum_{(\mu)}(-1)^{t+\eta_{w+t}}\kappa(\mu^{(1)})(x_{\sigma^{-1}(1)},\dots,\psi_{u-1}(\mu^{(2)};x_{\sigma^{-1}(l)},\dots),\dots)               \\
        -\sum_{(\mu)}(-1)^{w+\eta_1}\psi_{u-1}(\mu^{(1)};x_{\sigma^{-1}(1)},\dots,\kappa(\mu^{(2)})(x_{\sigma^{-1}(l)},\dots),\dots) \\
        -\sum_{i=1}^r(-1)^{w+\beta}\psi_u(\mu;x_1,\dots,d_1(x_i),\dots,x_r)
    \end{multline}
    for $0\leq u\leq w+1$, using the convention $\psi_{-1}=\psi_{w+1}=0$. The constants $\beta$ and $\eta$ are in \eqref{constants_beta_eta_alpha}.
    Moreover, the induced cohomological splitting
    \[\D{H}^n_{\O{O}}(A,M)=\prod_{w+t=n}\D{H}^{w,t}_{\O{O}}(A,M)\]
    is independent of the choice of horizontal resolution and is natural with respect to graded $\O{O}$-algebra morphisms and module morphisms.
\end{proposition}

\begin{proof}
    From the very definition of operadic cochain complex we have that
    \begin{align*}
        C_{\O{O}}^n(\bires{A},M)=\prod_{u\geq 0}\prod_{s\geq 0}\prod_{p\in\Z}\hom((\K{\O{O}})^{(u)}(\bires{A})_{s,p-s}, M_{p-n}).
    \end{align*}
    Here we use the weight grading of $\K{\O{O}}$ ($u\geq 0$) and the bigrading of $(\K{\O{O}})^{(u)}(\bires{A})$ in the statement, which is concentrated in non-negative horizontal degrees ($s\geq 0$). If we make the change of variables $s=w-u$ then $u\leq w$ since $s\geq 0$, hence
    \begin{align*}
        C_{\O{O}}^n(\bires{A},M)=\prod_{w\geq 0}\prod_{u=0}^w\prod_{p\in\Z}\hom((\K{\O{O}})^{(u)}(\bires{A})_{w-u,p-w+u}, M_{p-n}).
    \end{align*}
    Here, the middle product is actually direct a sum since it has finitely many factors.
    Calling $t=n-w$, we have that
    \[C^n_{\O{O}}(\bires{A},M)=\prod_{w+t=n}C^{w,t}_{\O{O}}(\bires{A},M).\]
    Note that $C^{w,t}_{\O{O}}(\bires{A},M)=0$ for $w<0$ because the direct sum is empty in this case.

    We now check that the operadic complex differential (co)restricts to
    \[d\colon C^{w,t}_{\O{O}}(\bires{A},M)\To C^{w+1,t}_{\O{O}}(\bires{A},M).\]
    This is a consequence of the following observations about formula \eqref{operadic_cohomology_differential_0}:
    \begin{itemize}
        \item The first term vanishes because $M$ has trivial differential.
        \item In the second term, the weight of $\mu^{(2)}$ is one less than the weight of $\mu$ because $\kappa$ vanishes in weight $\neq 1$.
        \item The same happens in the third term with $\mu^{(1)}$.
        \item In the last term, $d$ is the differential of $\bires{A}$, which is just $d=d_1$ by minimality. Hence, $d$ reduces the horizontal degree of $x_s$ by one.
    \end{itemize}
    This proves the chain complex splitting of $C^*_{\O{O}}(\bires{A},M)$ in the statement and the formula \eqref{operadic_cohomology_differential} for the differential of $C^{*,t}_{\O{O}}(\bires{A},M)$.

    The splitting is well defined in derived cohomology because different splittings associated to different horizontal resolutions compare well to the preferred horizontal resolution. Here we use the $E^2$-equivalence in Proposition \ref{horizontal_resolutions_exist}. Moreover, given a graded $\O{O}$-algebra morphism $f\colon B\to A$ the induced map in derived operadic cohomology $f^*\colon \D{H}^n_{\O{O}}(A,M)\to \D{H}^n_{\O{O}}(B,M)$ can be computed by using a compatible map between horizontal resolutions $g\colon\bires{B}\to\bires{A}$ as in Proposition \ref{horizontal_resolutions_exist}. Since this map is bigraded, $f^*$ preserves the splittings.
\end{proof}

\begin{remark}
    We will use several times the following instances of \eqref{operadic_cohomology_differential}. Let $\psi\in C^{w,t}_{\O{O}}(\bires{A},M)$. Recall that $(\K{\O{O}})^{(u)}$ is $\unit$, $sE$, and $s^2R$ for $u=0,1,2$, respectively. First,
    \begin{equation}\label{operadic_cohomology_differential_1}
        d(\psi)_u(1;x)=(-1)^{w+1}\psi_u(1;d_1(x)).
    \end{equation}
    Given $\mu\in E(r)$,
    \begin{multline}\label{operadic_cohomology_differential_E}
        d(\psi)_u(s\mu;x_1,\dots,x_r)
        =\sum_{i=1}^r(-1)^{t+\left(\beta+1\right)(w+t)}\mu(x_1,\dots,\psi_{u-1}(1;x_i),\dots,x_r)\\
        -(-1)^{w}\psi_{u-1}(1;\mu(x_1,\dots,x_r))
        +\sum_{i=1}^r(-1)^{w+\beta}\psi_u(s\mu;x_1,\dots,d_1(x_i),\dots,x_r),
    \end{multline}
    where $\beta$ is as in \eqref{constants_beta_eta_alpha}
    Given
    \begin{equation}\label{relation}
        \Gamma=\sum(\mu^{(1)}\circ_l\mu^{(2)})\cdot\sigma\in R(r)\subset (E\circ_{(1)}E)(r),
    \end{equation}
    \begin{multline}\label{operadic_cohomology_differential_R}
        d(\psi)_u(s^2\Gamma;x_1,\dots,x_r) =
        \sum(-1)^{\theta_{w,t}}\mu^{(1)}(x_{\sigma^{-1}(1)},\dots,\psi_{u-1}(s\mu^{(2)};x_{\sigma^{-1}(l)},\dots),\dots)               \\
        +\sum(-1)^{w+\delta}\psi_{u-1}(s\mu^{(1)};x_{\sigma^{-1}(1)},\dots,\mu^{(2)}(x_{\sigma^{-1}(l)},\dots),\dots) \\
        -\sum_{i=1}^r(-1)^{w+\gamma}\psi_u(s^2\Gamma;x_1,\dots,d_1(x_i),\dots,x_r).
    \end{multline}
    Here,
    \begin{equation}\label{constants}
        \begin{split}
            &\gamma         =\abs{\Gamma}+\sum_{j=1}^{i-1}\abs{x_j},                                \quad
            \delta         =\alpha_\sigma+\abs{\mu^{(2)}}\sum_{i=1}^{l-1}\abs{x_{\sigma^{-1}(i)}},\\
            &\theta_{w,t} =\alpha_\sigma+w+(\abs{\mu^{(1)}})(w+t+1)+(\abs{\mu^{(2)}}+1+w+t)\sum_{i=1}^{l-1}\abs{x_{\sigma^{-1}(i)}},
        \end{split}
    \end{equation}
    with $\alpha_\sigma$ as in \eqref{constants_beta_eta_alpha}. We have used that the infinitesimal decomposition of $\K{\O{O}}$ satisfies the following formulas,
    \begin{equation}\label{infinitesimal_decomposition}
        \begin{split}
            \Delta_{(1)}(1)         & =1\circ_11,                                           \\
            \Delta_{(1)}(s\mu)      & =1\circ_1(s\mu)+\sum_{i=1}^r(s\mu)\circ_i1,           \\
            \Delta_{(1)}(s^2\Gamma) & = 1\circ_1(s^2\Gamma)+\sum_{i=1}^r(s^2\Gamma)\circ_i1
            + \sum(-1)^{\abs{\mu^{(1)}}}((s\mu^{(1)})\circ_l(s\mu^{(2)}))\cdot\sigma.
        \end{split}
    \end{equation}
\end{remark}

\begin{example}\label{examples_cohomology}
    The bigraded derived operadic cohomology of algebras and modules with trivial differential can be computed as follows for the classical operads, compare \cite[12.4.1]{loday_vallette_2012_algebraic_operads}.
    \begin{enumerate}
        \item For differential graded associative algebras, operadic cohomology in positive degrees coincides with Hochschild's up to a shift. Moreover,
              \[C^{w,t}_{\O{A}}(\bires{A},M)=\bigoplus_{u=0}^w\prod_{\substack{p\in\Z\\i_0+\cdots+i_u+u=w\\j_0+\cdots+j_u=p-w}}\hom((\bires{A})_{i_0,j_0}\otimes\cdots\otimes (\bires{A})_{i_u,j_u}, M_{p-w-t})\]
              and the differential of $\psi=(\psi_0,\dots,\psi_w)\in C^{w,t}_{\O{A}}(\bires{A},M)$ is given by the following formula for $0\leq u\leq w+1$,
              \begin{align*}
                  d(\psi)_u(x_0,\dots, x_u)={} & (-1)^{w+(w+t+u-1)\abs{x_0}}x_0\psi_{u-1}(x_1,\dots, x_u)                               \\
                                               & -\sum_{i=0}^{u-1}(-1)^{w+i}\psi_{u-1}(x_0,\dots, x_ix_{i+1},\dots,x_u)                 \\
                                               & -(-1)^{w+u}\psi_{u-1}(x_0,\dots, x_{u-1})x_u                                           \\
                                               & -\sum_{i=0}^u(-1)^{w+u+\sum_{j=0}^{s-1}\abs{x_j}}\psi_u(x_0,\dots,d_1(x_i),\dots,x_u).
              \end{align*}
              Here we understand that $\psi_{-1}=0=\psi_{w+1}$.

        \item In the differential graded commutative case, we have Harrison's instead of Hochschild cohomology and
              \[C^{w,t}_{\O{C}}(\bires{A},M)\subset C^{w,t}_{\O{A}}(\bires{A},M)\]
              is the subcomplex of cochains $\psi=(\psi_0,\dots,\psi_w)$ such that each $\psi_u$ vanishes on shuffles in the sense of \cite[Example 5.18 (2)]{derived_homotopy_algebras}.
        \item In the case of differential graded Lie algebras, we obtain the Chevalley--Eilenberg cohomology. More precisely,
              \[C^{w,t}_{\O{C}}(\bires{A},M)\subset\bigoplus_{u=0}^w\prod_{\substack{p\in\Z\\i_0+\cdots+i_u+u=w\\j_0+\cdots+j_u=p-w}}\hom((\bires{A})_{i_0,j_0}\otimes\cdots\otimes (\bires{A})_{i_u,j_u}, M_{p-w-t})\]
              consists of the elements $\psi=(\psi_0,\dots,\psi_w)$ such that each $\psi_u$ is skew-symmetric in the sense of \cite[Example 5.18 (3)]{derived_homotopy_algebras} and, given $0\leq u\leq w+1$,
              \begin{multline*}
                  d(\psi)_u(x_0,\dots, x_u)=-\sum_{i=1}^u(-1)^{w+i+\abs{x_i}\sum\limits_{k=i+1}^u\abs{x_k}}[\psi_{u-1}(\dots,\widehat{x}_i,\dots),x_i]\\
                  -\sum_{0\leq i<j\leq u}(-1)^{w+i+j-1+(\abs{x_i}+\abs{x_j})\sum\limits_{k=0}^{i-1}\abs{x_k}+\abs{x_j}\sum\limits_{k=i+1}^{j-1}\abs{x_k}}\psi_{u-1}([x_i,x_j],\dots, \widehat{x}_i,\dots, \widehat{x}_j,\dots) \\
                  -\sum_{i=0}^u(-1)^{w+u+\sum_{j=0}^{s-1}\abs{x_j}}\psi_u(x_0,\dots,d_1(x_i),\dots,x_u)
              \end{multline*}
              with the convention $\psi_{-1}=0=\psi_{w+1}$.
    \end{enumerate}
\end{example}

\begin{remark}\label{error_sagave}
    Sagave's definition of derived cohomology for $\O{O}=\O{A}$ the associative operad does not coincide with ours.
    Roughly speaking, the reason is that he uses horizontal resolutions as coefficients, rather than graded modules.

    More precisely, Sagave uses a Hochschild-like cohomology $HH^{q,t}(E,F)$ for bicomplex $\O{A}$-algebras $E$ with coefficients in $E$-bimodules $F$ also in bicomplexes. Both $E$ and $F$ are required to have trivial vertical differential. He claims that his definition is invariant in both variables under $E^2$-equivalences $E\to F$ between $\kk$-projective minimal bicomplex $\O{A}$-algebras. Unfortunately this is not true, as we now show with an example.

    Sagave's $HH^{0,0}(E,E)$ coincides with the central elements of $E$ in bidegree $(0,0)$. The (horizontal) differential of $E$ does not play any role in horizontal degree $0$. Hence, we can take the unital bicomplex $\O{A}$-algebras $E=\Z$ concentrated in bidegree $(0,0)$ and $F=\Z[x,y]/(y^2)$ with $x$ of bidegree $(0,0)$, $y$ of bidegree $(1,0)$, and $d_1(y)=x$. Clearly, the unital map $E\to F$ is an $E^2$-equivalence. However $HH^{0,0}(E,E)=\Z$ but $HH^{0,0}(F,F)=\Z[x]$. Sagave's cohomology is not invariant in higher horizontal degrees either, however computing explicit examples becomes more complicated and it is not really worth to include them here.
\end{remark}

\section{Derived homotopy algebras}\label{section:derived_homotopy_algebras}

In this section we briefly recall from \cite{derived_homotopy_algebras} the theory of derived homotopy algebras. They are filtered homotopy algebras which allow the construction of minimal models over rings, even for algebras with torsion homology. This theory was initiated by Sagave \cite{sagave_2010_dgalgebras_derived_algebras} in the associative case, which attracted some attention  \cite{livernet_roitzheim_whitehouse_2013_derived_infty_algebras,cirici_egas_santander_livernet_whitehouse_2018_derived_infinity_algebras}, and then continued in \cite{maes_2016_derived_homotopy_algebras,derived_homotopy_algebras} over other operads.

Recall that a \emph{homotopy $\O{O}$-algebra} is an algebra over the operad $\Oinf{\O{O}}=\Omega\K{\O{O}}$, the cobar construction of the Koszul dual cooperad $\K{\O{O}}$ of $\O{O}$. 


\begin{remark}\label{homotopy_algebra}
    The structure of a homotopy $\O{O}$-algebra can be described in terms of the Koszul dual cooperad $\K{\O{O}}$.
    An $\Oinf{\O{O}}$-algebra is a complex $A$ equipped with structure morphisms
    \begin{equation*}
        \begin{split}
            \K{\O{O}}(r)_{n_0}\otimes A_{n_1}\otimes\cdots\otimes A_{n_r} & \To A_{n_0-1+n_1+\cdots+n_r},  \\
            \mu\otimes x_1\otimes\cdots\otimes x_r                                  & \;\mapsto\;\mu(x_1,\dots,x_r),
        \end{split}
    \end{equation*}
    satisfying $1(x)=0$ and certain formulas related to the infinitesimal decomposition of $\K{\O{O}}$, see \cite[\S10.1.2]{loday_vallette_2012_algebraic_operads} or \cite[Remark 5.2]{derived_homotopy_algebras}.
    %
    %
\end{remark}

Derived homotopy $\O{O}$-algebras are split filtered homotopy $\O{O}$-algebras. Let us make this structure precise.

\begin{definition}\label{derived_homotopy_algebra}
    A \emph{derived $\Oinf{\O{O}}$-algebra} or \emph{derived homotopy $\O{O}$-algebra} $A$ is an $\Oinf{\O{O}}$-algebra such that:
    \begin{enumerate}
        \item For each $n\in\Z$, the degree $n$ module splits as
              \[A_n=\bigoplus_{p+q=n}A_{p,q}\]
              with $A_{p,q}=0$ if $p<0$, like bicomplexes, see \S\ref{section:splitting_cohomology}.

        \item If we denote
              \[F_mA_n=\bigoplus_{\substack{p+q=n\\p\leq m}}A_{p,q},\]
              the differential of $A$ satisfies $d(F_mA_n)\subset F_mA_{n-1}$, so $A$ is a \emph{filtered complex}.

        \item The $\Oinf{\O{O}}$-algebra structure is compatible with the filtration, i.e.~the structure maps (co)restrict to
              \[\K{\O{O}}(r)_{n_0}\otimes F_{m_1}A_{n_1}\otimes\cdots\otimes F_{m_r}A_{n_r}\To F_{m_1+\cdots+m_r}A_{n_0-1+n_1+\cdots+n_r}.\]
    \end{enumerate}
\end{definition}

\begin{remark}\label{remark_dha}
    This is not the original definition in \cite[Definition 5.1]{derived_homotopy_algebras}, but an equivalent one from \cite[Corollary 5.21]{derived_homotopy_algebras}. The differential and the structure maps of $A$ are therefore determined by pieces, $i\geq 0$,
    \begin{align*}
        d_i\colon A_{p,q}                                                   & \To A_{p-i,q-1+i};                             \\
        \K{\O{O}}(r)_{s}\otimes A_{p_1,q_1}\otimes\cdots\otimes A_{p_r,q_r} & \To A_{p_1+\cdots+p_r-i,s-1+q_1+\cdots+q_r+i}, \\
        \mu\otimes x_1\otimes\cdots\otimes x_r                              & \;\mapsto\;\mu_i(x_1,\dots,x_r).
    \end{align*}
    The derived homotopy $\O{O}$-algebra equations can be rewritten in terms of these pieces, see \cite[Proposition 5.10]{derived_homotopy_algebras}. One of them is \begin{equation}\label{1_equation}
        1_i(x)=0,\qquad i\geq 0,
    \end{equation}
    where $1\in \kk=\unit(1)_0=(\K{\O{O}})^{(0)}(1)_0\subset\K{\O{O}}(1)_0$. 
    See Remark \ref{minimal_derived_homotopy_algebra_equations} below for some of these equations.
\end{remark}

\begin{remark}\label{classical_derived_homotopy_algebras}
    For the three classical operads (associative, commutative, and Lie), derived homotopy algebras look similar, since they have operations of the same bidegrees, but they satisfy quite different laws, see \cite[Example 5.18]{derived_homotopy_algebras} for full details.

    A derived homotopy associative algebra $A$ has operations $d_i$ and $m_{i,r}\colon A^{\otimes^r}\to A$ of bidegree $(-i,r-2+i)$ for $r\geq 2$ and $i\geq 0$. A derived homotopy commutative algebra too, and the $m_{i,r}$ vanish on shuffles. In the case of derived homotopy Lie algebras, these operations are skew-symmetric and denoted by $\ell_{i,r}\colon A^{\otimes^r}\to A$.

    If $\O{O}=\unit$ is the initial operad, an $\O{O}$-algebra is just a complex and a derived homotopy $\O{O}$-algebra is known as a \emph{twisted complex}. Its only structure maps are the $d_i$, $i\geq0$.
\end{remark}

\begin{remark}\label{infinity_morphisms}
    Usual morphisms between $\Oinf{\O{O}}$-algebras are not enough for many purposes, hence this category is often enlarged by means of the so-called $\infty$-morphisms.
    An \emph{$\infty$-morphism} between $\Oinf{\O{O}}$-algebras $f\colon A\leadsto B$ is given by structure morphisms
    \begin{align*}
        \K{\O{O}}(r)_{n_0}\otimes A_{n_1}\otimes \cdots\otimes A_{n_r} & \To B_{n_0+n_1+\cdots+n_r},       \\
        \mu\otimes x_1\otimes\cdots\otimes x_r                         & \;\mapsto\;f(\mu)(x_1,\dots,x_r),
    \end{align*}
    satisfying some equations connected with the decomposition and the infinitesimal decomposition of $\K{\O{O}}$, see \cite[\S10.2.3]{loday_vallette_2012_algebraic_operads} and \cite[Remark 5.23]{derived_homotopy_algebras}.


    The \emph{underlying morphism} of an $\infty$-morphism $f\colon A\leadsto B$ is the chain map $f(1)\colon A\to B$, where $1\in\K{\O{O}}$ is the element explained in Remark \ref{homotopy_algebra}.
\end{remark}

\begin{definition}
    Given two derived $\Oinf{\O{O}}$-algebras $A,B$ a \emph{derived $\infty$-morphism} $f\colon A\leadsto B$ is an $\infty$-morphism whose structure maps (co)restrict to
    \[\K{\O{O}}(r)_{n_0}\otimes F_{m_1}A_{n_1}\otimes\cdots\otimes F_{m_r}A_{n_r}\To F_{m_1+\cdots+m_r}A_{n_0+n_1+\cdots+n_r}.\]
\end{definition}

\begin{remark}\label{remark_dhm}
    The original definition is in \cite[Definition 5.24]{derived_homotopy_algebras} and this equivalent description corresponds to \cite[Corollary 5.27]{derived_homotopy_algebras}. The structure maps of $f\colon A\leadsto B$ are determined by pieces, $i\geq 0$,
    \begin{align*}
        \K{\O{O}}(r)_s\otimes A_{p_1,q_1}\otimes\cdots\otimes A_{p_r,q_r} & \To B_{p_1+\cdots+p_r-i,s+q_1+\cdots+q_r+i}, \\
        \mu\otimes x_1\otimes\cdots\otimes x_r                            & \;\mapsto\; f(\mu)_i(x_1,\dots,x_r).
    \end{align*}
    The derived $\infty$-morphism equations can be rewritten in terms of these pieces, see \cite[Proposition 5.25]{derived_homotopy_algebras}. Some of these equations will be recalled in Remark \ref{special_derived_homotopy_algebra_equations} below.
\end{remark}

\begin{remark}\label{bicomplex_algebras_as_derived_homotopy_algebras}
    A bicomplex $\O{O}$-algebra is exactly the same thing as a derived homotopy $\O{O}$-algebra such that $d_i=0$ for $i\geq 2$ and
    \[\mu_i(x_1,\dots,x_r)=0\]
    if $i\geq 1$ or if $\mu\in\K{\O{O}}$ has weight $\neq 1$. Given $\mu\in E(r)$, the formula
    \begin{equation}\label{derived_infinity_or_not}
        (s\mu)_0(x_1,\dots,x_r)=\mu(x_1,\dots,x_r)
    \end{equation}
    relates the derived homotopy $\O{O}$-algebra structure, on the left, and the bicomplex $\O{O}$-algebra structure, on the right. 

    A differential graded $\O{O}$-algebra is the same as a bicomplex $\O{O}$-algebra concentrated in horizontal degree $0$, so $d_1=0$ and $d=d_0$.

    A bicomplex $\O{O}$-algebra morphism is the same as a derived $\infty$-morphism $f$ such that $f(\mu)_i=0$ for $i\geq 1$ or $\mu$ of weight $\geq 1$, i.e.~the only possibly non-trivial part is $f(1)_0$.
    %
\end{remark}

\begin{definition}\label{frito_variado}
    A derived $\infty$-morphism between derived homotopy $\O{O}$-algebras $f\colon A\leadsto B$ is an \emph{$E^2$-equivalence} if the filtration-preserving chain map $f(1)_0\colon A\to B$ induces an isomorphism between the $E^2$ pages of the spectral sequences associated to the filtrations.

    A derived homotopy $\O{O}$-algebra is \emph{minimal} if $d_0=0$. If $A$ is a differential graded $\O{O}$-algebra, a \emph{(derived) minimal model} is an $E^2$-equivalence $\minimal{A}\leadsto A$ whose source $\minimal{A}$ is a minimal and $\kk$-projective derived homotopy $\O{O}$-algebra.
\end{definition}

Under our standing assumptions, minimal models exist, and there is actually a preferred one which allows for comparison.

\begin{theorem}[{\cite[Theorem 7.25 and Proposition 7.26]{derived_homotopy_algebras}}]\label{compare_minimal_models}
    Any differential graded $\O{O}$-algebra $A$ has a minimal model $\minimal{A}'\leadsto A$ such that, if $\minimal{A}\leadsto A$ is another minimal model, then there exists an $E^2$-equivalence $\minimal{A}\leadsto \minimal{A}'$ which induces the identity in $H_*(A)$ on the $E^2$ terms of the corresponding spectral sequences. Moreover, if $f\colon B\to A$ is an $\O{O}$-algebra morphism and $\minimal{B}\leadsto B$ is a minimal model then there exists a derived $\infty$-morphism $\minimal{B}\leadsto \minimal{A}'$ which induces $H_*(f)$ on $E^2$ terms.
\end{theorem}

In \cite[Examples 6.6, 6.8, and 6.11]{derived_homotopy_algebras} we explicitly compute minimal models for some of the algebras in Examples \ref{dugger_shipley}, \ref{commutative}, and \ref{lie} below. 


\begin{remark}\label{remark_minimal}\label{E2_term_minimal_model}
    In a minimal model, $d_1$ is a bidegree $(-1,0)$ differential on $\minimal{A}$ by minimality. We call it \emph{horizontal differential}. Its homology is the $E^2$ term of the spectral sequence of $\minimal{A}$, so it is isomorphic to $H_*(A)$ concentrated in horizontal degree $0$ since $f$ is an $E^2$-equivalence (in particular the complexes defined by $d_1$ are exact in positive horizontal degrees). This isomorphism is induced by $f(1)_{0}\colon \minimal{A}\to A$, which is trivial in positive horizontal degrees for degree reasons and maps horizontal degree $0$ elements in $\minimal{A}$ to cycles in $A$ by minimality.
\end{remark}

\begin{remark}\label{underlying_bicomplex_algebra_of_minimal_derived_homotopy_algebra}
    By \cite[Remark 6.3]{derived_homotopy_algebras}, any minimal derived homotopy $\O{O}$-algebra has an underlying minimal bicomplex $\O{O}$-algebra structure given by restricting to $d_0=0$, $d_1$, and the operations $(s\mu)_0$, $\mu\in E$, i.e.~\eqref{derived_infinity_or_not}, see Remark \ref{bicomplex_algebras_as_derived_homotopy_algebras}. Minimality is essential here. Moreover, if $f\colon A\leadsto B$ is a derived $\infty$-morphism then $f(1)_0\colon A\to B$ is a morphism between the underlying bicomplex $\O{O}$-algebras.
\end{remark}

\begin{remark}\label{horizontal_from_minimal}
    Given a differential graded $\O{O}$-algebra $A$, any minimal model $f\colon\minimal{A}\leadsto A$ produces a horizontal resolution $\rho\colon \bires{A}\to H_*(A)$ of the homology graded $\O{O}$-algebra of $A$. Here, $\bires{A}$ is the bicomplex $\O{O}$-algebra underlying $\minimal{A}$ in the sense of Remark \ref{underlying_bicomplex_algebra_of_minimal_derived_homotopy_algebra} and $\rho$ is the composite
    \[\bires{A}\to Z_*(A)\twoheadrightarrow H_*(A)\]
    where $Z_*(A)\subset A$ are the cycles, the second arrow is the natural projection onto homology, and the first arrow is the corestriction of $f(1)_0\colon \minimal{A}\to A$, see Remark \ref{remark_minimal}.

    In particular, if $f\colon\minimal{A}\leadsto A$ and $f'\colon\minimal{A}'\leadsto A$ are minimal models and $g\colon \minimal{A}\leadsto\minimal{A}'$ is an $E^2$-equivalence inducing the identity in the $E^2$ terms of the corresponding spectral sequences (it is $H_*(A)$ concentrated in horizontal degree $0$ in both cases), then the induced horizontal resolutions satisfy
    \[\rho'g(1)_0=\rho.\]
\end{remark}

\begin{remark}\label{minimal_derived_homotopy_algebra_equations}
    The following are some of the equations of a \emph{minimal} derived homotopy $\O{O}$-algebra in terms of the operations in Remark \ref{remark_dha}. These are the ones we will later need. We will use the underlying bicomplex $\O{O}$-algebra structure in Remark \ref{underlying_bicomplex_algebra_of_minimal_derived_homotopy_algebra}, and more specifically the notation \eqref{derived_infinity_or_not} for $\mu\in E$.

    By minimality $d_0=0$ and $d_1^2=0$. We also have,
    \begin{equation}\label{twisted_complex_equation_3}
        d_1d_2+d_2d_1=0.
    \end{equation}
    Recall that $(\K{\O{O}})^{(u)}$ is $sE,s^2R$ for $u=1,2$, respectively.
    Given $\mu\in E(r)$,
    \begin{multline}\label{big_equation_E_2}
        0=d_1((s\mu)_1(x_1,\dots,x_r))+d_2(\mu(x_1,\dots,x_r))\\
        -\sum_{i=1}^r(-1)^{\beta}\left(\mu(x_1,\dots,d_2(x_i),\dots,x_r)
        +(s\mu)_1(x_1,\dots,d_1(x_i),\dots,x_r)\right)
    \end{multline}
    where $\beta$ is as in \eqref{constants_beta_eta_alpha}.
    Given
    $\Gamma\in R(r)$ like in \eqref{relation},
    we have
    \begin{multline}\label{big_equation_R_1}
        0=d_1((s^2\Gamma)_0(x_1,\dots,x_r))
        +\sum_{i=1}^r(-1)^{\gamma}(s^2\Gamma)_0(x_1,\dots,d_1(x_i),\dots,x_r)\\
        -\sum_{(\mu)}(-1)^{\delta}\left(\mu^{(1)}(x_{\sigma^{-1}(1)},\dots,(s\mu^{(2)})_1(x_{\sigma^{-1}(l)},\dots),\dots)\right.\\
        \left.+(s\mu^{(1)})_1(x_{\sigma^{-1}(1)},\dots,\mu^{(2)}(x_{\sigma^{-1}(l)},\dots),\dots)\right)
    \end{multline}
    where $\gamma$ and $\delta$ been defined in \eqref{constants}
    These formulas use \cite[Proposition 5.10]{derived_homotopy_algebras}, \eqref{1_equation}, and the fact that the infinitesimal decomposition of $\K{\O{O}}$ satisfies \eqref{infinitesimal_decomposition}.
    With the Sweedler notation in \eqref{sweedler_infinitesimal_decomposition}, we also have
    \begin{equation}\label{big_equation_O_0}
        0=\sum_{(\mu)}(-1)^{\eta_1}\mu^{(1)}_0(x_{\sigma^{-1}(1)},\dots,\mu^{(2)}_0(x_{\sigma^{-1}(l)},\dots),\dots),
    \end{equation}
    where $\eta_1$ is as in \eqref{constants_beta_eta_alpha}. 
\end{remark}

\begin{remark}\label{special_derived_homotopy_algebra_equations}
    We will need the following equations of a derived $\infty$-morphism $f\colon A\leadsto B$ between derived homotopy $\O{O}$-algebras $A, B$ in terms of the operations in Remark \ref{remark_dhm}. Below $A$ will always be minimal, and $B$ will either be also minimal or it will be a differential graded $\O{O}$-algebra,  regarded as a derived homotopy $\O{O}$-algebra concentrated in horizontal degree $0$, see Remark \ref{bicomplex_algebras_as_derived_homotopy_algebras}. Here we also use the underlying bicomplex $\O{O}$-algebra structure of a minimal derived homotopy $\O{O}$-algebra, see  Remark \ref{underlying_bicomplex_algebra_of_minimal_derived_homotopy_algebra} and the notation \eqref{derived_infinity_or_not} for $\mu\in E$.

    Recall that $(\K{\O{O}})^{(u)}$ is $\unit,sE,s^2R$ for $u=0,1,2$, respectively.

    Let $B$ be a differential graded $\O{O}$-algebra. Then
    \begin{align}
        \label{derived_infinity_morphism_1_0}
        df(1)_0                              & =0, \\
        \label{derived_infinity_morphism_1_1}
        df(1)_1(x)-f(1)_0d_1(x)              & =0, \\
        \label{derived_infinity_morphism_1_2}
        df(1)_2(x)-f(1)_0d_2(x)-f(1)_1d_1(x) & =0.
    \end{align}
    Given $\mu\in E(r)$,
    \begin{multline}\label{derived_infinity_morphism_E_1}
        df(s\mu)_1(x_1,\dots,x_r)
        +\sum_{i=1}^r(-1)^{\beta}f(s\mu)_0(x_1,\dots,d_1(x_i),\dots,x_r)
        \\
        =f(1)_0(s\mu)_1(x_1,\dots,x_r)+f(1)_1\mu(x_1,\dots,x_r)\\\
        -\sum_{i=1}^r\mu(f(1)_0(x_1),\dots,f(1)_1(x_i),\dots,f(1)_0(x_r)),
    \end{multline}
    where $\beta$ has been defined in \eqref{constants_beta_eta_alpha}.

    Suppose now that $B$ is a minimal derived homotopy $\O{O}$-algebra, like $A$. Let us rename $f\colon A\leadsto B$ as $g\colon A\leadsto B$ because we will later use the formulas below with this name for the derived $\infty$-morphism. We have
    \begin{equation}\label{derived_infinity_morphism_eq_1_2}
        d_1g(1)_1(x)+d_2g(1)_0(x)-g(1)_0d_2(x)-g(1)_1d_1(x)=0.
    \end{equation}
    Moreover, for $\mu\in E(r)$,
    \begin{multline}\label{derived_infinity_morphism_eq_E_1}
        d_1g(s\mu)_0(x_1,\dots,x_r)+\sum_{i=1}^r(-1)^\beta g(s\mu)_0(x_1,\dots,d_1(x_i),\dots,x_r)=\\
        g(1)_0((s\mu)_1(x_1,\dots,x_r))+g(1)_1(\mu(x_1,\dots,x_r))\\
        -\sum_{i=1}^r\mu(g(1)_0(x_1),\dots,g(1)_1(x_i),\dots,g(1)_0(x_r))-(s\mu)_1(g(1)_0(x_1),\dots,g(1)_0(x_r))
    \end{multline}
    For $\Gamma\in R(r)$ like in \eqref{relation},
    \begin{multline}\label{derived_infinity_morphism_eq_R_0}
        0=g(1)_0((s^2\Gamma)_0(x_1,\dots,x_r))-(s^2\Gamma)_0(g(1)_0(x_1),\dots,g(1)_0(x_r))\\
        -\sum(-1)^{\delta}g(s\mu^{(1)})_0(x_{\sigma^{-1}(1)},\dots,\mu^{(2)}(x_{\sigma^{-1}(l)},\dots),\dots)\\
        -\sum(-1)^{\varepsilon}\mu^{(1)}(g(1)_0(x_{\sigma^{-1}(1)}),\dots,g(s\mu^{(2)})_0(x_{\sigma^{-1}(l)},\dots),\dots),
    \end{multline}
    where
    \begin{equation}\label{epsilon}
        \varepsilon = \alpha_\sigma+\abs{\mu^{(1)}}+(\abs{\mu^{(2)}}+1)\sum_{i=1}^{l-1}\abs{x_{\sigma^{-1}(i)}},
    \end{equation}
    and $\alpha_\sigma$ and $\delta$ have been defined in \eqref{constants_beta_eta_alpha} and \eqref{constants}, respectively.

    All this follows from \cite[Proposition 5.25]{derived_homotopy_algebras}, \eqref{1_equation}, and the fact that the infinitesimal decomposition $\Delta_{(1)}\colon\K{\O{O}}\to \K{\O{O}}\circ_{(1)}\K{\O{O}}$ and the decomposition $\Delta\colon\K{\O{O}}\to \K{\O{O}}\circ\K{\O{O}}$ of the Koszul dual cooperad $\K{\O{O}}$ satisfy \eqref{infinitesimal_decomposition} and
    \begin{align}
        \nonumber
        \Delta(1)         ={} & (1;1),                              \\
        \nonumber
        \Delta(s\mu)      ={} & (1;s\mu)+(s\mu;1,\dots,1),          \\
        \nonumber
        \Delta(s^2\Gamma) ={} & (1;s^2\Gamma)+(s^2\Gamma;1,\dots,1)
        \nonumber
        +\sum(-1)^{\abs{\mu^{(1)}}}(s\mu^{(1)};1,\stackrel{l-1}{\dots},1,s\mu^{(2)},1,\dots,1)\cdot\sigma.
    \end{align}
\end{remark}

\section{Derived universal Massey products}\label{section:universal_massey_products}

We are finally ready to define the derived cohomology class this paper takes its name from.


\begin{definition}\label{universal_massey_product}
    Given an $\O{O}$-algebra $A$ in $\chain$, its \emph{derived universal Massey product}
    \[\g{A}{\O{O}}=\{\m\}\in \D{H}^{2,-1}_{\O{O}}(H_*(A),H_*(A))\]
    is the derived operadic cohomology class represented by the cochain $\m$ defined in the following way. Choose a minimal model $f\colon \minimal{A}\leadsto A$, consider the induced horizontal resolution $\rho\colon \bires{A}\to H_*(A)$ in Remark \ref{horizontal_from_minimal}, and define the representing cocycle $\m=(\m_0,\m_1,\m_2)\in C_{\O{O}}^{2,-1}(\bires{A},H_*(A))$
    by the following formulas,
    \begin{align*}
        \m_0(1;x)                     & = \rho(d_2(x)),                                            \\
        \m_1(s\mu;x_1,\dots,x_r)      & = \rho((s\mu)_1(x_1,\dots,x_r)),\quad \mu\in E(r),         \\
        \m_2(s^2\Gamma;x_1,\dots,x_r) & = \rho((s^2\Gamma)_0(x_1,\dots,x_r)),\quad \Gamma\in R(r).
    \end{align*}
    See Proposition \ref{splitting} for notation and recall that $(\K{\O{O}})^{(u)}$ is $\unit$, $sE$, and $s^2R$ for $u=0,1,2$, respectively.
\end{definition}


\begin{remark}\label{classical_derived_universal_massey_product}
    With the notation in Example \ref{examples_cohomology} and Remark \ref{classical_derived_homotopy_algebras}, the derived universal Massey product of a differential graded associative algebra $A$ with minimal model $f\colon \minimal{A}\leadsto A$ is represented by the bidegree $(2,-1)$ cochain $(\rho d_2,\rho m_{1,2},\rho m_{0,3})$. The same formula holds in the commutative case, and also in the Lie case replacing $m$ with $\ell$. The case $\O{O}=\unit$ will be considered in \S\ref{section:chain_complexes} below.

    Notice that the associative derived universal Massey product is indeed defined by the same pieces of a minimal model as in \cite[Proposition 5.4]{sagave_2010_dgalgebras_derived_algebras}, but composed with $\rho$ and living in a different cohomology theory, see Remark \ref{error_sagave}.
\end{remark}


The following result shows that the derived universal Massey product is a well-defined natural invariant of differential graded $\O{O}$-algebras.

\begin{theorem}\label{dimitrova_like}
    In the setting of Definition \ref{universal_massey_product}, the following statements hold:
    \begin{enumerate}
        \item The cochain $\m$ in is indeed a cocycle.
        \item Its cohomology class is independent of the choice of minimal model.
        \item Given an $\O{O}$-algebra morphism $f\colon A\to B$, the morphisms induced by $H_*(f)$ in derived operadic cohomology
              \[\D{H}^{2,-1}_{\O{O}}(H_*(A),H_*(A)){\To}\D{H}^{2,-1}_{\O{O}}(H_*(A),H_*(B)){\longleftarrow}\D{H}^{2,-1}_{\O{O}}(H_*(B),H_*(B))\]
              take $\g{A}{\O{O}}$ and $\g{B}{\O{O}}$ to the same class in the middle.
        \item Given a weight-preserving morphism of operads satisfying our standing assumptions $f\colon \O{P}\to\O{O}$, the morphism induced by $f$ on derived operadic cohomology
              \[H_{\O{O}}^{2,-1}(H_*(A),H_*(A))\To H_{\O{P}}^{2,-1}(H_*(A),H_*(A))\]
              takes $\g{A}{\O{O}}$ to $\g{A}{\O{P}}$.
    \end{enumerate}
\end{theorem}

\begin{proof}
    We start with (1). We will see that $d(\m)=0$ as a consequence of the minimal  derived homotopy $\O{O}$-algebra equations for $\minimal{A}$. By \eqref{operadic_cohomology_differential_1},
    \[
        d(\m)_0(1;x)={}  -\m_0(1;d_1(x))
        ={}               -\rho (d_2d_1(x))
        ={}               \rho d_1d_2(x)
        ={}               0.
    \]
    Here we also use \eqref{twisted_complex_equation_3} and the fact that $\rho d_1=0$  (see Remark \ref{certainty}).

    Given $\mu\in E(r)$, using \eqref{operadic_cohomology_differential_E} and the definition of $\m$
    \begin{multline*}
        d(\m)_1(s\mu;x_1,\dots,x_r)
        ={}                                  \sum_{i=1}^r(-1)^{\beta}\underbrace{\mu(x_1,\dots, \rho (d_2(x_i)),\dots,x_r)}_{\rho (\mu(x_1,\dots, d_2(x_i),\dots,x_r))}\\
        -\rho d_2(\mu(x_1,\dots,x_r))
        +\sum_{i=1}^r(-1)^{\beta}\rho ((s\mu)_1(x_1,\dots,d_1(x_i),\dots,x_r)).
    \end{multline*}
    Here, in the underbrace we use that $H_*(A)$ is regarded as an $\bires{A}$-module via the bicomplex $\O{O}$-algebra morphism $\rho \colon\bires{A}\to H_*(A)$. This vanishes by \eqref{big_equation_E_2} since $\rho d_1=0$.

    For $\Gamma\in R(r)$ as in \eqref{relation}, using now \eqref{operadic_cohomology_differential_R}
    \begin{multline*}
        d(\m)_2(s^2\Gamma;x_1,\dots,x_r)={} \\
        \sum(-1)^{\theta_{2,-1}}\underbrace{\mu^{(1)}(x_{\sigma^{-1}(1)},\dots,\rho ((s\mu^{(2)})_1(x_{\sigma^{-1}(l)},\dots)),\dots)}_{\rho (\mu^{(1)}(x_{\sigma^{-1}(1)},\dots,(s\mu^{(2)})_1(x_{\sigma^{-1}(l)},\dots),\dots))}                                                         \\
        +\sum(-1)^{\delta}\rho ((s\mu^{(1)})_1(x_{\sigma^{-1}(1)},\dots,\mu^{(2)}(x_{\sigma^{-1}(l)},\dots),\dots))              \\
        -\sum_{i=1}^r(-1)^{\gamma}\rho ((s^2\Gamma)_0(x_1,\dots,d_1(x_i),\dots,x_r)).
    \end{multline*}
    This vanishes by \eqref{big_equation_R_1} since $\rho d_1=0$ and $\theta_{2,-1}=\delta$.

    Finally, for $\mu\in (\K{\O{O}})^{(3)}$, using \eqref{operadic_cohomology_differential} and the notation \eqref{sweedler_infinitesimal_decomposition},
    \begin{align*}
        d(\m)_3(\mu;x_1,\dots,x_r) ={} &
        \sum_{(\mu)}(-1)^{\eta_{1}-1}\kappa(\mu^{(1)})(x_{\sigma^{-1}(1)},\dots,\m_2(\mu^{(2)};x_{\sigma^{-1}(l)},\dots),\dots)                                                                                                     \\
                                       & - \sum_{(\mu)}(-1)^{\eta_1}\m_2(\mu^{(1)};x_{\sigma^{-1}(1)},\dots,\kappa(\mu^{(2)})(x_{\sigma^{-1}(l)},\dots),\dots)                                                                      \\
        ={}                            &
        -\sum_{(\mu)'}(-1)^{\eta_1}\underbrace{\mu^{(1)}_0(x_{\sigma^{-1}(1)},\dots,\rho (\mu^{(2)}_0(x_{\sigma^{-1}(l)},\dots)),\dots)}_{\rho (\mu^{(1)}_0(x_{\sigma^{-1}(1)},\dots,\mu^{(2)}_0(x_{\sigma^{-1}(l)},\dots),\dots))} \\
                                       & - \sum_{(\mu)''}(-1)^{\eta_1}\rho (\mu^{(1)}_0(x_{\sigma^{-1}(1)},\dots,\mu^{(2)}_0(x_{\sigma^{-1}(l)},\dots),\dots))                                                                      \\
        ={}                            &
        -\sum_{(\mu)}(-1)^{\eta_1}\rho \mu^{(1)}_0(x_{\sigma^{-1}(1)},\dots,\mu^{(2)}_0(x_{\sigma^{-1}(l)},\dots),\dots).
    \end{align*}
    Here, $(\mu)'$ (resp.~$(\mu)''$) runs over the summands of \eqref{sweedler_infinitesimal_decomposition} with $\mu^{(1)}$ (resp.~$\mu^{(2)}$) of weight $1$. We use that the weights of $\mu^{(1)}$ and $\mu^{(2)}$ add up to $3$ and \eqref{1_equation}. This vanishes by \eqref{big_equation_O_0}.

    Let us tackle (2). It suffices to show that, in the situation of Theorem \ref{compare_minimal_models}, the cocycles $\m$ and $\m'$ defined by $f\colon\minimal{A}\leadsto A$ and $f'\colon\minimal{A}'\leadsto A$ give rise to the same cohomology class. In order to compare them, we use the $E^2$-equivalence $g\colon \minimal{A}\leadsto\minimal{A}'$ in that theorem. The underlying morphism of bicomplex $\O{O}$-algebras $g(1)_0\colon\bires{A}\to\bires{A}'$ induces a map of complexes,
    \[g(1)_0^*\colon C_{\O{O}}^{*,t}(\bires{A}',H_*(A))\To C_{\O{O}}^{*,t}(\bires{A},H_*(A))\]
    which takes $\m'$ to
    \begin{align*}
        g(1)_0^*(\m')_0(1;x)                     & = \rho '(d_2(g(1)_0(x))),                                                    \\
        g(1)_0^*(\m')_1(s\mu;x_1,\dots,x_r)      & = \rho '((s\mu)_1(g(1)_0(x_1),\dots,g(1)_0(x_r))),\quad \mu\in E(r),         \\
        g(1)_0^*(\m')_2(s^2\Gamma;x_1,\dots,x_r) & = \rho '((s^2\Gamma)_0(g(1)_0(x_1),\dots,g(1)_0(x_r))),\quad \Gamma\in R(r).
    \end{align*}
    We define the cochain $\psi\in C_{\O{O}}^{1,-1}(\bires{A},H_*(A))$ as
    \begin{align*}
        \psi_0(1;x)                & = \rho '(g(1)_1(x)),                               \\
        \psi_1(s\mu;x_1,\dots,x_r) & = \rho '(g(s\mu)_0(x_1,\dots,x_r)),\quad \mu\in E.
    \end{align*}
    Let us check that $d(\psi)=g(1)_0^*(\m')-\m$.

    By \eqref{operadic_cohomology_differential_1} and \eqref{derived_infinity_morphism_eq_1_2},
    \begin{align*}
        d(\psi)_0(1;x) 
         & =  \rho '(g(1)_1d_1(x))                           \\
         & =  \rho '(d_1g(1)_1(x)+d_2g(1)_0(x)-g(1)_0d_2(x)) \\
         & = \rho 'd_2g(1)_0(x)-\rho d_2(x)                  \\
         & = g(1)_0^*(\m')_0(1;x)-\m_0(1;x)
    \end{align*}
    Here we use that $\rho 'd_1=0$ like above. We also use that $\rho 'g(1)_0=\rho $, see Remark \ref{horizontal_from_minimal}.

    Given $\mu\in E(r)$, by \eqref{operadic_cohomology_differential_E} and \eqref{derived_infinity_morphism_eq_E_1},
    \begin{multline*}
        d(\psi)_1(s\mu;x_1,\dots,x_r)=
        -\sum_{i=1}^r\underbrace{\mu(x_1,\dots, \rho '(g(1)_1(x_i)),\dots,x_r)}_{\rho '(\mu(g(1)_0(x_1),\dots,g(1)_1(x_i),\dots,g(1)_0(x_r)))}\\
        +\rho '(g(1)_1(\mu(x_1,\dots,x_r)))
        -\sum_{i=1}^r(-1)^{\beta}\rho '(g(s\mu)_0(x_1,\dots,d_1(x_i),\dots,x_r))                                                                \\
        =\rho 'd_1g(s\mu)_0(x_1,\dots,x_r)-\rho '(g(1)_0(s\mu)_1(x_1,\dots,x_r)\\
        +(s\mu)_1(g(1)_0(x_1),\dots,g(1)_0(x_r)))\\
        =-\rho ((s\mu)_1(x_1,\dots,x_r))+\rho '((s\mu)_1(g(1)_0(x_1),\dots,g(1)_0(x_r)))\\
        =-\m_1(s\mu;x_1,\dots,x_r)+g(1)_0^*(\m')_1(s\mu;x_1,\dots,x_r).
    \end{multline*}
    Here we use again the equations $\rho 'd_1=0$ and $\rho 'g(1)_0=\rho $, and the fact that $\rho '$ is a bicomplex $\O{O}$-algebra morphism.

    Given $\Gamma\in R(r)$, by \eqref{operadic_cohomology_differential_R} and \eqref{derived_infinity_morphism_eq_R_0},
    \begin{multline*}
        d(\psi)_2(s^2\Gamma;x_1,\dots,x_r)={}\\
        \sum(-1)^{\theta_{1,-1}}\underbrace{\mu^{(1)}(x_{\sigma^{-1}(1)},\dots,\rho '(g(s\mu^{(2)})_0(x_{\sigma^{-1}(l)},\dots)),\dots)}_{\rho '(\mu^{(1)}(g(1)_0(x_{\sigma^{-1}(1)}),\dots,g(s\mu^{(2)})_0(x_{\sigma^{-1}(l)},\dots),\dots))}                                                 \\
        +\sum(-1)^{1+\delta}\rho '(g(s\mu)_0(x_{\sigma^{-1}(1)},\dots,(s\mu^{(2)})_0(x_{\sigma^{-1}(l)},\dots),\dots))\\
        = -\rho '(g(1)_0((s^2\Gamma)_0(x_1,\dots,x_r))+(s^2\Gamma)_0(g(1)_0(x_1),\dots,g(1)_0(x_r)))\\
        = -\rho ((s^2\Gamma)_0(x_1,\dots,x_r))+\rho '(s^2\Gamma)_0(g(1)_0(x_1),\dots,g(1)_0(x_r))\\
        = -\m_2(s^2\Gamma;x_1,\dots,x_r)+g(1)_0^*(\m')_2(s^2\Gamma;x_1,\dots,x_r).
    \end{multline*}
    Here we use once again $\rho 'd_1=0$, $\rho 'g(1)_0=\rho $, and that $\rho '$ is a bicomplex $\O{O}$-algebra morphism. We also use that $\theta_{1,-1}=\epsilon+1$, see \eqref{constants} and \eqref{epsilon}.

    Part (3) is a slight generalization of (2) which also uses Theorem \ref{compare_minimal_models}, concretely the last part. Hence we skip it.

    Lastly, in (4) any minimal model for $A$ as an $\O{O}$-algebra can be restricted to a minimal model as a $\O{P}$-algebra and the statement holds already at the level of the cochains defined by a given minimal model. Such a restriction is possible because the operad map $f$ preserves weights, so it gives rise to a cooperad map between their Koszul duals $\K{\O{P}}\to\K{\O{O}}$.
\end{proof}

The derived universal Massey product is an obstruction to formality.

\begin{corollary}
    If $A$ is a formal differential graded $\O{O}$-algebra then $\g{A}{\O{O}}=0$.
\end{corollary}

\begin{proof}
    By Theorem \ref{dimitrova_like}, we can directly assume that $A=H_*(A)$, i.e.~$A$ has trivial differential. Then a horizontal resolution is also a minimal model for $A$, compare Remark \ref{bicomplex_algebras_as_derived_homotopy_algebras}. The cocycle $\m$ in Definition \ref{universal_massey_product} is trivial if $\minimal{A}$ is just a bicomplex $\O{O}$-algebra, hence we are done.
\end{proof}

\begin{remark}
    If $\kk$ is a field, deriving operadic cohomology is unnecessary, see Corollary \ref{corolario_todo}. Moreover, we can construct a minimal model for any differential graded $\O{O}$-algebra $A$ on its homology $H_*(A)$, and the derived universal Massey product coincides with Dimitrova's class since representing cocycles agree on the nose, compare \cite[\S4]{muro_2021_massey_products_algebras}. The same happens for $\kk$ a commutative ring if $H_*(A)$ is $\kk$-projective, compare \cite[Remark 4.2]{muro_2021_massey_products_algebras}.
\end{remark}

\section{Chain complexes}\label{section:chain_complexes}

The case of $\O{O}=\unit$ the initial operad does not lack of interest. Its algebras are chain complexes, its derived homotopy algebras are twisted complexes, and its derived $\infty$-morphisms are twisted maps. We do not need any hypotheses on the ground ring because this operad is non-symmetric.

We here compute the derived cohomology of a graded module and the derived universal Massey product of a chain complex. Any operad $\O{O}$ satisfying our assumptions admits an operad map $\unit\to\O{O}$, the unit. By naturality, it is sometimes possible to check that the derived universal Massey product of an $\O{O}$-algebra is non-trivial by just looking at its underlying complex, see Theorem \ref{dimitrova_like}. We will see below, in Example \ref{example_sagave}, that this is the case for the differential graded associative algebra in \cite[Example 5.7]{sagave_2010_dgalgebras_derived_algebras}.

A graded $\unit$-algebra $A$ is just a graded module and a horizontal resolution $\bires{A}$ is a sequence of projective resolutions of the modules $A_n$, $n\in\Z$. Moreover, a graded $A$-module $M$ is just another graded module, which can be completely unrelated to $A$ for there is not structure intertwining both. The factors of the operadic cochain complex in Proposition \ref{splitting} are
\[C_{\unit}^{w,t}(\bires{A},M)=\prod_{q\in\Z}\hom((\bires{A})_{w,q},M_{q-t})\]
with differential $\partial\colon C_{\unit}^{w,t}(E,M)\to C_{\unit}^{w+1,t}(E,M)$ is $d(\psi)(x)=-(-1)^wfd_1(x)$.
Therefore, the derived operadic cohomology of $A$ with coefficients in $M$ is
\[\D{H}_{\unit}^{w,t}(A,M)=\prod_{q\in\Z}\ext_{\kk}^w(A_q,M_{q-t}).\]
In particular, if $A$ is now a chain complex (regarded as a $\unit$-algebra), its derived universal Massey product lies in
\[\D{H}_{\unit}^{2,-1}(H_*(A),H_*(A))=\prod_{q\in\Z}\ext_{\kk}^2(H_q(A),H_{q+1}(A)).\]

We have the following intrinsic characterization of the derived universal Massey product of a chain complex.

\begin{proposition}\label{extensions}
    Given a chain complex $A$,
    \[\cdots\to A_{q+1}\stackrel{d}{\To}A_q\to\cdots,\]
    its derived universal Massey product $\g{A}{\unit}$ is represented by the opposite of the following sequence of extensions, $q\in\Z$,
    \[H_{q+1}(A)\hookrightarrow\frac{A_{q+1}}{d(A_{q+2})}\stackrel{d}{\To}\ker[A_{q}\stackrel{d}{\to} A_{q-1}]\twoheadrightarrow H_q(A).\]
\end{proposition}

\begin{proof}
    Let $f\colon\minimal{A}\leadsto A$ be a minimal model, which consists of maps, $i\geq 0$,
    \[f_i=f(1)_i\colon (\minimal{A})_{i,n-i}\To A_n\]
    satisfying certain equations. The derived universal Massey product is the sequence of extension classes represented by the following composites, $q\in\Z$,
    \[(\minimal{A})_{2,q}\stackrel{d_2}{\To}(\minimal{A})_{0,q+1}\stackrel{f_0}{\To}Z_{q+1}(A)\twoheadrightarrow H_{q+1}(A),\]
    see Definition \ref{universal_massey_product} and Remark \ref{horizontal_from_minimal}.

    Let us consider the following piece of the twisted map $f\colon \minimal{A}\leadsto A$,
    \begin{center}
        \begin{tikzcd}
            A_{q+2} \ar[d, "d"'] &&&\\
            A_{q+1} \ar[d, "d"'] & (\minimal{A})_{0,q+1} \ar[l, "f_0"']
            & (\minimal{A})_{1,q+1} \ar[l, "d_1"', near start] 
            & \\
            A_{q} \ar[d, "d"'] & (\minimal{A})_{0,q} \ar[l, "f_0"']
            & (\minimal{A})_{1,q} \ar[l, "d_1"] \ar[llu, "f_1"']
            & (\minimal{A})_{2,q} \ar[l, "d_1"] \ar[llu, "d_2", near end] 
            \ar[llluu, "f_2"', near end]\\
            A_{q-1} & (\minimal{A})_{0,q-1} \ar[l, "f_0"] & (\minimal{A})_{1,q-1} \ar[l, "d_1"] \ar[llu, "f_1", near end] & (\minimal{A})_{2,q-1} \ar[l, "d_1"] \ar[llu, "d_2", near end] \ar[llluu, "f_2", very near end]
        \end{tikzcd}
    \end{center}
    The following diagram commutes by \eqref{derived_infinity_morphism_1_0}, \eqref{derived_infinity_morphism_1_1}, and \eqref{derived_infinity_morphism_1_2},
    \begin{center}
        \begin{tikzcd}
            (\minimal{A})_{2,q} \ar[r,"d_1"] \ar[d,"d_2"] & (\minimal{A})_{1,2} \ar[r, "d_1"] \ar[ddd, "f_1"] & (\minimal{A})_{0,q} \ar[ddd, "f_0"] \ar[r, two heads] & H_q(A) \ar[ddd, equal] \\
            (\minimal{A})_{0,q+1} \ar[d,"-f_0"] &&& \\
            Z_{q+1}(A) \ar[d,two heads] &&& \\
            H_{q+1}(A) \ar[r, hook] & \frac{A_{q+1}}{d(A_{q+2})} \ar[r, "d"] & \ker [d\colon A_q\to A_{q-1}] \ar[r, two heads] & H_q(A)
        \end{tikzcd}
    \end{center}
    Here, on top, we find part of the horizontal resolution $\bires{A}$ of $H_q(A)$ arising from the minimal model $\minimal{A}$ as in Remark \ref{horizontal_from_minimal}. The bottom line is the extension in the statement. This proves the proposition.
\end{proof}


\begin{example}\label{example_sagave}
    The example considered by Sagave in \cite[Example 5.7]{sagave_2010_dgalgebras_derived_algebras} over $\kk=\Z/(p^2)$ is the unital differential graded associative algebra $A=\kk[x]/(x^2)$ with $\abs{x}=1$ and $d(x)=p$. The underlying chain complex is
    \[\cdots\to 0\to \Z/(p^2)\stackrel{p}{\To}\Z/(p^2)\to 0\to\cdots\]
    concentrated in degrees $1$ and $2$. By Proposition \ref{extensions}, $\g{A}{\unit}$ is represented by the extension,
    \[\Z/(p)\hookrightarrow \Z/(p^2)\stackrel{p}{\To}\Z/(p^2)\twoheadrightarrow\Z/(p).\]
    It is well known that this extension generates
    \[\D{H}_{\unit}^{2,-1}(H_*(A),H_*(A))=\ext_{\Z/(p^2)}^2(\Z/(p),\Z/(p))\cong \Z/(p).\]
    Hence $\g{A}{\unit}\neq0$ and $0\neq \g{A}{\O{A}}\in\D{H}^{2,-1}_{\O{A}}(H_*(A),H_*(A))$ by the naturality of derived universal Massey products with respect to restriction of scalars along the unit $\unit\to\O{A}$.

    If we replace $\Z$ and $p\in\Z$ with $\Q[t]$ and $t\in \Q[t]$ then $A$ is a differential graded commutative algebra with $0\neq \g{A}{\O{O}}\in \D{H}^{2,-1}_{\O{O}}(H_*(A),H_*(A))$ for $\O{O}=\unit,\O{A},\O{C}$.
\end{example}

\section{Torsion Massey products}\label{section:torsion_massey}

In this section we give an easy criterion to show that some $\O{O}$-algebras over rings like $\kk=\Z,\Q[t]$ have non-trivial derived universal Massey product. The theory in \S\ref{section:chain_complexes} does not apply here because these rings have global dimension $\leq 1$.

We will consider the following secondary operation in the cohomology of a differential graded $\O{O}$-algebra.

\begin{definition}\label{torsion_massey_product}
    Let $\mu\in E(r)$ be an arity $r$ generator of our operad $\O{O}$ and $t_i\in \kk$, $1\leq i\leq r$, scalars such that $\sum_{i=1}^rt_i=0$. Given a differential graded $\O{O}$-algebra $A$ and $x_i\in H_{*}(A)$ with $t_ix_i=0$, $1\leq i\leq r$, choose representatives $y_i\in A_{\abs{x_i}}$ of the $x_i$ and elements $z_i\in A_{\abs{x_i}+1}$ with
    \[d_i(z_i)=t_iy_i,\]
    for all $1\leq i\leq r$. The \textit{torsion Massey product}
    \begin{equation}\label{stp}
        \secop{\mu}{t_1,\dots,t_r}{x_1,\dots,x_r}\in \frac{H_{\abs{\mu}+\abs{x_1}+\cdots+\abs{x_r}+1}(A)}{\sum_{i=1}^r\mu(x_1,\dots, x_{i-1},H_{\abs{x_i}+1}(A),x_{i+1},\dots,x_r)}
    \end{equation}
    is the element represented by the homology class of
    \begin{equation}\label{representative}
        \sum_{i=1}^r(-1)^{\beta}\mu(y_1,\dots, y_{i-1},z_i,y_{i+1},\dots,y_r)\in A_{\abs{\mu}+\abs{x_1}+\cdots+\abs{x_r}+1},
    \end{equation}
    where $\beta$ is as in \eqref{constants_beta_eta_alpha}.
    The denominator in \eqref{stp} is known as the \emph{indeterminacy} of the torsion Massey product.
\end{definition}

\begin{remark}
    The torsion Massey product is the operadic Massey product, in the sense of \cite{muro_2021_massey_products_algebras}, associated to the following relation in $\O{O}$,
    \begin{equation*}
        \sum_{i=1}^r\mu\circ_i(t_i\cdot 1)=\sum_{i=1}^rt_i(\mu\circ_i 1)=\left(\sum_{i=1}^r t_i\right)\mu=0
    \end{equation*}
    Strictly speaking, this is not a relation in $\O{O}$ in the sense we are using this word in this paper, i.e.~it is not given by an element of $R$. However, this kind of relation suffices to define an operadic Massey product, see \cite[Remark 4.2]{muro_2021_massey_products_algebras}. In particular, the different choices of $y_i,z_i$ produce all possible elements of the coset
    \begin{equation*}
        \secop{\mu}{t_1,\dots,t_r}{x_1,\dots,x_r}\subset H_{\abs{\mu}+\abs{x_1}+\cdots+\abs{x_r}+1}(A),
    \end{equation*}
    Moreover, torsion Massey products are invariant under differential graded $\O{O}$-algebra morphisms $f\colon A\to B$, in the sense that the induced map $f_*\colon H_*(A)\to H_*(B)$ satisfies
    \begin{equation*}
        f(\secop{\mu}{t_1,\dots,t_r}{x_1,\dots,x_r})\subset \secop{\mu}{t_1,\dots,t_r}{f(x_1),\dots,f(x_r)}\subset H_{\abs{\mu}+\abs{x_1}+\cdots+\abs{x_r}+1}(B)
    \end{equation*}
    whenever the left hand side is defined. As a consequence, if $A$ is formal all torsion Massey products contain $0$.


    The condition on the $t_i$ implies that $t_r=-\sum_{i=1}^{r-1}t_i$, hence we can omit $t_r$ from notation,
    \[\secop{\mu}{t_1,\dots,t_{r-1}}{x_1,\dots,x_r}=\secop{\mu}{t_1,\dots,t_{r-1},t_r}{x_1,\dots,x_r}.\]
    In the examples below we will consider cases where the generator $\mu$ has arity $r=2$. In this case $\secop{\mu}{t}{x_1,x_2}$ is defined whenever $t\cdot x_1=0=t\cdot x_2$ for $t\in\kk$.

    Torsion Massey products can only be non-trivial for $r\geq 2$. If $t_i=0$ for all $i$ then the corresponding torsion Massey products vanish because then the homology class of \eqref{representative} lies in the indeterminacy.

\end{remark}

Now we show how to compute torsion Massey products from a minimal model. Actually, from the pieces of a minimal model defining the derived universal Massey product.

\begin{theorem}\label{torsion_massey_infinity}
    Suppose that we have a differential graded $\O{O}$-algebra $A$, a minimal model $f\colon \minimal{A}\leadsto A$, a generator $\mu\in E(r)$, and elements $x_i\in H_*(A)$ and $t_i\in\kk$ such that $t_ix_i=0$, $1\leq i\leq r$, and $\sum_{i=1}^rt_i=0$. Choose elements $u_i\in (\minimal{A})_{0,\abs{x_i}}$ such that $f(1)_0(u_i)\in A_{\abs{x_i}}$ represents $x_i$, $v_i\in (\minimal{A})_{1,\abs{x_i}}$ satisfying $d_1(v_i)=t_iu_i$, $1\leq i\leq r$, and
    \[w\in (\minimal{A})_{2,\abs{\mu}+\sum_{i=1}^r\abs{x_i}}\]
    such that
    \[d_1(w)=\sum_{i=1}^r(-1)^{\beta}(s\mu)_0(u_1,\dots,u_{i-1},v_i,u_{i+1},\dots,u_r),\]
    where $\beta$ is as in \eqref{constants_beta_eta_alpha}. Then, for $\rho$ as in Remark \ref{horizontal_from_minimal},
    \[\sum_{i=1}^r(-1)^{\beta}\rho(s\mu)_1(u_1,\dots,v_i,\dots,u_r)-\rho d_2(w)\in H_{\abs{\mu}+\sum_{i=1}^r\abs{x_i}+1}(A)\]
    belongs to the torsion Massey product $\secop{\mu}{t_1,\dots,t_r}{x_1\dots,x_r}$.
\end{theorem}

\begin{proof}
    First, notice that the elements $u_i$, $v_i$, and $w$ must exist (although they are not unique) by Remark \ref{remark_minimal}. Since $f(1)_0$ is a chain map, we can choose
    \[y_i=f(1)_0(u_i)\]
    for the computation of the torsion Massey product following Definition \ref{torsion_massey_product}.

    By \eqref{derived_infinity_morphism_1_1},
    \[df(1)_1(v_i)=t_if(1)_0(u_i),\]
    so we can choose
    \[z_i=f(1)_1(v_i).\]

    By \eqref{derived_infinity_morphism_E_1} for $(x_1,\dots,x_r)=(u_1,\dots,u_{i-1},v_i,u_{i+1},\dots,u_r)$, since $d_1(u_j)=0$ and $f(1)_1(u_j)=0$ for degree reasons,
    \begin{multline*}
        df(s\mu)_1(u_1,\dots,v_i,\dots,u_r)
        +(-1)^{\beta}t_if(s\mu)_0(u_1,\dots,u_r)\\
        =f(1)_0(s\mu)_1(u_1,\dots,v_i,\dots,u_r)+f(1)_1(s\mu)_0(u_1,\dots,v_i,\dots,u_r)\\
        -\mu(y_1,\dots,z_i,\dots,y_r).
    \end{multline*}

    With the previous choices, using this last formula, the cycle in \eqref{representative} is
    \begin{multline*}
        \sum_{i=1}^r(-1)^{\beta}\left(f(1)_0(s\mu)_1(u_1,\dots,v_i,\dots,u_r)+f(1)_1(s\mu)_0(u_1,\dots,v_i,\dots,u_r)\right)\\
        -\sum_{i=1}^rt_if(s\mu)_0(u_1,\dots,u_r)-\sum_{i=1}^r(-1)^{\beta}df(s\mu)_1(u_1,\dots,v_i,\dots,u_r).
    \end{multline*}
    The first summation in the second line vanishes because $\sum_{i=1}^rt_i=0$. Moreover, the last summation in the second line is a boundary, so if we discard it we obtain another cycle in the same homology class. Hence, the torsion Massey product contains the homology class of
    \begin{equation*}
        \sum_{i=1}^r(-1)^{\beta}f(1)_0(s\mu)_1(u_1,\dots,v_i,\dots,u_r)
        +f(1)_1d_1(w).
    \end{equation*}

    By \eqref{derived_infinity_morphism_1_2}, this cycle coincides with
    \begin{equation*}
        \sum_{i=1}^r(-1)^{\beta}f(1)_0(s\mu)_1(u_1,\dots,v_i,\dots,u_r)
        +df(1)_2(w)-f(1)_0d_2(w).
    \end{equation*}
    We can discard the middle term because it is a boundary. The homology class of this cycle is the one in the statement, see Remark \ref{horizontal_from_minimal}
\end{proof}

\begin{remark}\label{massey_infinity}
    In \cite[\S2]{muro_2021_massey_products_algebras}, given a differential graded $\O{O}$-algebra $A$, a relation $\Gamma\in R(r)$ in $\O{O}$, and $x_1,\dots,x_r\in H_*(A)$ satisfying certain vanishing conditions we defined a coset
    \[\langle x_1,\dots,x_r\rangle_\Gamma\subset H_{\abs{\Gamma}+\sum_{i=1}^r\abs{x_i}+1}(A)\]
    that we called \emph{Massey product}, since it generalizes Massey's and Retakh's in the associative and Lie cases, respectively.

    In \cite[\S3]{muro_2021_massey_products_algebras}, we established a connection between these Massey products and non-derived minimal models whenever they exist (under strong projectivity and cofibrancy assumptions). Derived minimal models always exist (as recalled in Theorem \ref{compare_minimal_models}) and essentially the same proof shows that, given elements $u_i\in (\minimal{A})_{0,\abs{x_i}}$ such that $f(1)_0(u_i)\in A_{\abs{x_i}}$ represents $x_i$,
    \[f(1)_0((s^2\Gamma)_0(u_1,\dots,u_r))\in A_{\abs{\Gamma}+\sum_{i=1}^r\abs{x_i}+1}\]
    is a cycle whose homology class belongs to $\langle x_1,\dots,x_r\rangle_\Gamma$. Hence, derived minimal models also compute these \emph{non-torsion Massey products}.
\end{remark}

\begin{remark}
    Theorem \ref{torsion_massey_infinity} and Remark \ref{massey_infinity} assert that torsion Massey products and operadic Massey products associated to elements of $R$ can be computed from the cocycle representing the derived universal Massey product in Definition \ref{universal_massey_product}. We can actually use any representative. More precisely, given a differential graded $\O{O}$-algebra $A$, a horizontal resolution $\rho\colon\bires{A}\to H_*(A)$ and a representative $\psi\in C^{2,-1}_{\O{O}}(\bires{A},H_*(A))$ of $\g{A}{\O{O}}$, we can replace $\minimal{A}$, $\rho(s\mu)_1$, $\rho d_2$, and $\rho(s^2\Gamma)_0$ in Theorem \ref{torsion_massey_infinity} and Remark \ref{massey_infinity} with $\bires{A}$, $\psi_1(s\mu;\dots)$, $\psi_0(1;\dots)$, and $\psi_2(s^2\Gamma;\dots)$, respectively. This is tedious but straightforward to check, compare \cite[Proposition 4.1]{muro_2021_massey_products_algebras}. We here include a full proof of the following weaker statement.
\end{remark}

\begin{proposition}
    If $A$ is a differential graded $\O{O}$-algebra with trivial derived universal Massey product, then all torsion Massey products vanish in $H_*(A)$.
\end{proposition}

\begin{proof}
    Let $\psi\in C_{\O{O}}^{1,-1}(\bires{A},A)$ be a cochain such that satisfying $d(\psi)=\varphi$. With the notation in Theorem \ref{torsion_massey_infinity}, using \eqref{operadic_cohomology_differential_E}, \eqref{operadic_cohomology_differential_1}, and that $d_1(u_i)=0$ and $\psi_0(1;u_i)=0$ for degree reasons,
    \begin{multline*}
        \sum_{i=1}^r(-1)^{\beta}d(\psi)_1(s\mu;u_1,\dots,v_i,\dots,u_r)-d(\psi)_0(1;w)=\\
        -\sum_{i=1}^r(-1)^{\beta}\mu(u_1,\dots,\psi_0(1;v_i),\dots,u_r)
        +\sum_{i=1}^r(-1)^{\beta}\psi_0(1;\mu(u_1,\dots,v_i,\dots,u_r))\\
        -\sum_{i=1}^r\underbrace{\psi_1(s\mu;u_1,\dots,d_1(v_i),\dots,u_r)}_{t_i\psi_1(s\mu;u_1,\dots,u_r)}-\psi_0(1;d_1(w))=\\
        -\sum_{i=1}^r(-1)^{\beta}\mu(x_1,\dots,\psi_0(1;v_i),\dots,x_r).
    \end{multline*}
    Here we use that $\sum_{i=1}^rt_i=0$, and the fact that $A$ is regarded as an $\bires{A}$-module via the bicomplex $\O{O}$-algebra map $f\colon\bires{A}\to A$, which takes $u_i$ to $x_i$. Clearly, this element belongs to the indeterminacy.
\end{proof}

In order to conclude this paper, we now exhibit several examples of differential graded $\O{O}$-algebras over rings of global dimension $1$, for $\O{O}$ the classic operads, with some some non-vanishing torsion Massey product, and hence non-trivial derived universal Massey product. In particular, these differential graded $\O{O}$-algebras are not formal. In these examples, we cannot detect the non-triviality of derived universal Massey products by using operadic Massey products associated to elements of $R$. These are either not defined except in trivial situations or they vanish for degree reasons.

\begin{example}\label{dugger_shipley}
    Let $\O{O}=\O{A}$ be the associative operad with the usual presentation, with generating $\SS$-module $E$ concentrated in arity $2$ and degree $0$, generated by $\mu$, the operation representing the associative product. Hence $E(2)$ is a free $\SS_2$-module of rank $1$. Moreover, $R\subset E\circ_{(1)}E$ is generated by the associativity relation $\mu\circ_1\mu-\mu\circ_2\mu$.

    It would be natural to consider torsion Massey products associated to $\mu\in E$, and so we do in Example \ref{commutative} below, but we here consider the commutator
    \[\ell=\mu-\mu\cdot (1\; 2)\in E.\]

    Assume $A$ is a differential graded associative algebra which in addition is unital. Suppose that there exists $t\in\kk$ such that the unit $1\in H_*(A)$ satisfies $t\cdot 1=0$. Hence $t\cdot H_*(A)$, so
    \[\secop{\ell}{t}{1,x}\in\frac{H_{\abs{x}+1}}{[H_1(A),x]}\]
    is defined for all $x\in H_*(A)$. Here
    \[[y,z]=yz-(-1)^{\abs{y}\abs{z}}zy\]
    denotes the commutator bracket, which is the operation corresponding to $\mu$. For the computation of the indeterminacy we have used that $1\in H_*(A)$ is in the graded center. If $H_*(A)$ is commutative the indeterminacy vanishes.

    For $p\in\Z=\kk$ a prime, Dugger and Shipley consider in \cite{dugger_shipley_2009_curious_example_triangulatedequivalent} the differential graded unital associative algebra
    \[A=\frac{\Z\langle e,x^{\pm1}\rangle}{(e^2,ex+xe-x^2)},\qquad \abs{e}=\abs{x}=1,\]
    with differential
    \begin{align*}
        d(e) & =p, & d(x) & =0.
    \end{align*}
    Its homology is
    \[H_*(A)=\Z/(p)\langle x^{\pm1}\rangle,\]
    which is commutative (always in the graded sense) if and only if $p=2$. We will not distinguish between cycles in $A$ and their homology classes so as not to overload notation.

    We are going to compute $\secop{\ell}{p}{1,x}$. For $p\neq 2$,
    \[[x,x]=2x^2\in H_2(A)\]
    is a generator so $H_2(A)=[H_1(A),x]$ and the torsion Massey product trivially vanishes. However, for $p=2$ the indeterminacy vanishes by commutativity and the second defining relation of $A$ shows that
    \[\secop{\ell}{2}{1,x}=x^2\in H_2(A).\]

    The algebra $A$ is actually not formal for any prime $p$. This follows from \cite{dugger_shipley_2009_curious_example_triangulatedequivalent}. An alternative proof based on derived minimal models is given in \cite[Example 8.8]{derived_homotopy_algebras}.
\end{example}

\begin{example}\label{commutative}
    Let $\kk=\Q[t]$ and let $\O{O}=\O{C}$ be the commutative operad with the usual presentation, which is like the presentation of $\O{A}$ in Example \ref{dugger_shipley} but requiring $\mu\cdot(1\;2)=\mu$, so $E(2)$ is a free $\kk$-module of rank $1$ with trivial action of $\SS_2$.

    We consider the graded commutative (non-unital) algebra $A$ in \cite[Example 6.8]{derived_homotopy_algebras} generated by \[x,y,x_t,y_t\] in degrees $\abs{x}=\abs{y}=2$ and $\abs{x_t}=\abs{y_t}=3$, subject to the relations \[x^2,y^2,xy,xx_t,yy_t,x_ty.\] We endow it with the differential defined by
    \begin{align*}
        d(x)   & =d(y)=0, &
        d(x_t) & =tx,     &
        d(y_t) & =ty.
    \end{align*}
    Its homology $H_*(A)$ has basis $\{x,y,xy_t\}$ over $\Q$, not over $\kk$. Again, we do not distinguish between cycles in $A$ and their classes. Actually, $t\cdot H_*(A)=0$, so torsion Massey products associated to $\mu$ and $t$ are always defined.  The product in $H_*(A)$ is trivial, hence the indeterminacy always vanishes. Nevertheless,
    \[\secop{\mu}{t}{x,y}=-xy_t\in H_5(A).\]

    The relations in the previous differential graded algebra were designed so as to have a fully-computed small homology algebra. This was necessary in \cite{derived_homotopy_algebras}. However, here we can work without relations. More precisely, if $B$ is the free graded commutative algebra with the same generators (and no relations) and we define the differential in the same way, then
    \[\secop{\mu}{t}{x,y}=x_ty-xy_t\in H_5(B)\cong\Q\]
    is a generator and the indeterminacy vanishes because $H_3(B)=0$.

    Using the canonical operad morphism $\O{A}\to\O{C}$ we see that these examples work in exactly the same way when regarded as associative algebras.
\end{example}

\begin{example}\label{lie}
    Consider again $\kk=\Q[t]$ and let $\O{O}=\O{L}$ be the Lie operad with the standard presentation,
    with $E$ concentrated in arity $2$ and degree $0$, generated by $\ell$, the operation representing the Lie bracket, which must satisfy $\ell\cdot(1\;2)=-\ell$. Hence $E(2)$ is a free $\kk$-module of rank $1$ equipped with the sign action of $\SS_2$. The $\SS$-module of relations $R\subset E\circ_{(1)}E$ is generated by the Jacobi relation $(\ell\circ_1\ell)\cdot[()+(1\;2\;3)+(3\;2\;1)]$.

    The algebra $A$ we now consider, introduced in \cite[Example 6.11]{derived_homotopy_algebras}, behaves similarly to that in Example \ref{commutative}, despite they are algebras over different operads. It is a graded Lie algebra generated by \[x,y,x_t,y_t\] in the same degrees as above subject to the relations saying that all triple brackets vanish, and also all binary brackets except for $[x,y_t],[x_t,y_t]$ and their symmetrics. Then $H_*(A)$ has $\Q$-linear basis $\{x,y,[x,y_t]\}$ and $t\cdot H_*(A)=0$, so torsion Massey products associated to $\mu$ and $t$ are always defined. The graded Lie algebra $H_*(A)$ is abelian, hence indeterminacies vanish, and
    \[\secop{\mu}{t}{x,y}=-[x,y_t]\in H_5(A).\]

    We can also work without relations in this case, since the necessity of full homology computations in \cite{derived_homotopy_algebras} is not present in this paper. More precisely, if $B$ is the free graded Lie algebra with the same generators as the previous $A$ and we define the differential in the same way, then
    \[\secop{\mu}{t}{x,y}=[x_t,y]-[x,y_t]\in H_5(B)\cong\Q\]
    is a generator and the indeterminacy vanishes because $H_3(B)=0$.
\end{example}








\end{document}